\documentclass[12pt]{article}
\usepackage{amsfonts}
\usepackage{amssymb}
\title{Vertex operator algebras
associated to type $B$ affine Lie algebras
on admissible half-integer levels}
\author{Ozren Per\v{s}e}
\date{}
\begin{document}
\def \Z{\Bbb Z}
\def \C{\Bbb C}
\def \R{\Bbb R}
\def \Q{\Bbb Q}
\def \N{\Bbb N}
\def \tr{{\rm tr}}
\def \span{{\rm span}}
\def \Res{{\rm Res}}
\def \End{{\rm End}}
\def \E{{\rm End}}
\def \Ind {{\rm Ind}}
\def \Irr {{\rm Irr}}
\def \Aut{{\rm Aut}}
\def \Hom{{\rm Hom}}
\def \mod{{\rm mod}}
\def \ann{{\rm Ann}}
\def \<{\langle}
\def \>{\rangle}
\def \t{\tau }
\def \a{\alpha }
\def \e{\epsilon }
\def \l{\lambda }
\def \L{\Lambda }
\def \g{\gamma}
\def \b{\beta }
\def \om{\omega }
\def \o{\omega }
\def \c{\chi}
\def \ch{\chi}
\def \cg{\chi_g}
\def \ag{\alpha_g}
\def \ah{\alpha_h}
\def \ph{\psi_h}
\def \be{\begin{equation}\label}
\def \ee{\end{equation}}
\def \bl{\begin{lem}\label}
\def \el{\end{lem}}
\def \bt{\begin{thm}\label}
\def \et{\end{thm}}
\def \bp{\begin{prop}\label}
\def \ep{\end{prop}}
\def \br{\begin{rem}\label}
\def \er{\end{rem}}
\def \bc{\begin{coro}\label}
\def \ec{\end{coro}}
\def \bd{\begin{de}\label}
\def \ed{\end{de}}
\def \pf{{\bf Proof. }}
\def \voa{{vertex operator algebra}}

\newtheorem{thm}{Theorem}[section]
\newtheorem{prop}[thm]{Proposition}
\newtheorem{coro}[thm]{Corollary}
\newtheorem{conj}[thm]{Conjecture}
\newtheorem{lem}[thm]{Lemma}
\newtheorem{rem}[thm]{Remark}
\newtheorem{de}[thm]{Definition}
\newtheorem{hy}[thm]{Hypothesis}
\makeatletter \@addtoreset{equation}{section}
\def\theequation{\thesection.\arabic{equation}}
\makeatother \makeatletter

\newcommand{\binom}[2]{{{#1}\choose {#2}}}
    \newcommand{\nno}{\nonumber}
    \newcommand{\lbar}{\bigg\vert}
    \newcommand{\p}{\partial}
    \newcommand{\dps}{\displaystyle}
    \newcommand{\bra}{\langle}
    \newcommand{\ket}{\rangle}
 \newcommand{\res}{\mbox{\rm Res}}
\renewcommand{\hom}{\mbox{\rm Hom}}
  \newcommand{\epf}{\hspace{2em}$\Box$}
 \newcommand{\epfv}{\hspace{1em}$\Box$\vspace{1em}}
\newcommand{\nord}{\mbox{\scriptsize ${\circ\atop\circ}$}}
\newcommand{\wt}{\mbox{\rm wt}\ }

\maketitle
\begin{abstract}
Let $L(n-l+\frac{1}{2},0)$ be the vertex operator algebra associated to
an affine Lie algebra of type $B_{l}^{(1)}$
at level $n-l+\frac{1}{2}$, for a positive integer~$n$.
We classify irreducible
$L(n-l+\frac{1}{2},0)$-modules and show that every
$L(n-l+\frac{1}{2},0)$-module is completely reducible. In the special
case $n=1$, we study a category of weak
$L(-l+\frac{3}{2},0)$-modules which are in the category $\cal{O}$
as modules for the associated affine Lie algebra. We classify
irreducible objects in that category and prove semisimplicity
of that category.
\end{abstract}

\footnotetext[2]{
Partially supported by the
Ministry of Science, Education and Sports
of the Republic of Croatia, grant 0037125.}

\section{Introduction}

Let ${\frak g}$ be a simple finite-dimensional Lie algebra and
$\hat{\frak g}$ the associated affine Lie algebra. For any complex
number $k$, denote by $L(k,0)$ the irreducible highest weight
$\hat{\frak g}$-module with the highest weight $k \Lambda_{0}$. Then
$L(k,0)$ has a natural vertex operator algebra structure for any
$k\ne -h^{\vee}$. The representation theory of $L(k,0)$ heavily
depends on the choice of level $k \in \C$. If $k$ is a positive integer,
$L(k,0)$ is a rational vertex operator algebra (cf. [FZ], [Z]),
i.e. the category of $\Z _{+}$-graded weak $L(k,0)$-modules is
semisimple. Irreducible objects in that category
are integrable highest weight $\hat{\frak g}$-modules
of level $k$ ([FZ],[L]). The corresponding
associative algebra $A(L(k,0))$, defined in [Z],
is finite-dimensional (cf. [KWn]).
In some cases such as $k \notin \Q$ or $k < -h^{\vee}$
(studied in [KL1] and [KL2]),
categories of $L(k,0)$-modules have significantly different
structure then categories of $L(k,0)$-modules for a positive
integer $k$. But there are examples of rational levels $k$ such that
the category of weak $L(k,0)$-modules which are in the category $\cal{O}$
as $\hat{\frak g}$-modules, has similar structure as the
category of $\Z _{+}$-graded weak $L(k,0)$-modules for positive integer levels $k$.
These are so called admissible levels, defined by V. Kac and M.
Wakimoto (cf. [KW1] and [KW2]). D. Adamovi\'{c} (cf. [A1] and
[A2]) studied vertex operator algebras associated to affine Lie
algebras of type $C_{l}^{(1)}$ on admissible half-integer levels.
D. Adamovi\'{c} and A. Milas [AM], and C. Dong, H.-S. Li and G.
Mason [DLM] studied vertex operator algebras associated to affine
Lie algebras of type $A_{1}^{(1)}$ on all admissible levels. It is
shown that in these cases vertex operator algebra $L(k,0)$ has finitely
many irreducible weak modules from the category $\cal{O}$ and that
every weak $L(k,0)$-module from the category $\cal{O}$ is
completely reducible. One can say that these vertex operator
algebras are rational in the category $\cal{O}$.
In [AM], authors gave a conjecture that
vertex operator algebras $L(k,0)$, for all admissible levels $k$, are
rational in the category $\cal{O}$.
In this paper we give examples of a vertex operator algebras $L(k,0)$
on admissible levels $k$ for which we prove some parts of the conjecture
from [AM]. Admissible modules for affine Lie algebras were
also recently studied in [A3], [FM], [GPW], [W].

We consider the case of an affine Lie algebra
of type $B_{l}^{(1)}$ and the corresponding vertex operator
algebra $L(n-l+\frac{1}{2},0)$, for any positive integer $n$. We show that
$n-l+\frac{1}{2}$ is an admissible level for this affine Lie algebra.
The results on admissible modules from [KW1] imply that $L(n-l+\frac{1}{2},0)$
is a quotient of the generalized Verma module
by the maximal ideal generated by a singular vector.
Using results from [Z], [FZ], we can identify the corresponding
associative algebra $A(L(n-l+\frac{1}{2},0))$ with a certain
quotient of $U({\frak g})$. Algebra $A(L(n-l+\frac{1}{2},0))$
is infinite-dimensional in this case. Using
methods from [MP], [A2] we get that irreducible
$A(L(n-l+\frac{1}{2},0))$-modules from the category $\mathcal{O}$ are
in one-to-one correspondence with zeros of the certain set of
polynomials ${\mathcal P}_{0}$. By calculating certain polynomials
from that set we obtain the classification of irreducible
finite-dimensional $A(L(n-l+\frac{1}{2},0))$-modules,
and in the special case $n=1$, the classification of irreducible
$A(L(-l+\frac{3}{2},0))$-modules from the category $\mathcal{O}$.
Using results from [Z], we obtain the classification of irreducible
$L(n-l+\frac{1}{2},0)$-modules, and in the special case $n=1$,
the classification of irreducible weak
$L(-l+\frac{3}{2},0)$-modules from the category $\mathcal{O}$.
Using these classifications and results from [KW2], we show
that every $L(n-l+\frac{1}{2},0)$-module is completely reducible,
and in the case $n=1$, that every weak
$L(-l+\frac{3}{2},0)$-module from the category $\mathcal{O}$
is completely reducible.

The method for classification of irreducible $L(k,0)$-modules
used in this paper depends on a relatively simple formula
for the singular vector in the generalized Verma module. For a
general admissible level $k$, a more global method for
classification is needed.

The author expresses his gratitude to Professors
Dra\v{z}en Adamovi\'{c} and Mirko Primc for their helpful
advice and constant support.

\section{Vertex operator algebras associated to affine Lie algebras}

This section is preliminary. We recall some necessary definitions
and fix the notation. We review certain results about vertex
operator algebras and corresponding modules. The emphasis is on
the class of vertex operator algebras associated to affine Lie
algebras, because we study a special case in that class in
Sections 3 and 4.

\subsection{Vertex operator algebras and modules}

Let $(V, Y, {\bf 1}, \omega)$ be a vertex operator algebra
(cf. [B], [FHL] and [FLM]).
An {\it ideal} in a vertex operator algebra $V$ is a subspace $I$
of $V$ satisfying $ Y(a,z)I \subseteq I[[z,z^{-1}]] $
for any $a \in V$. Given an ideal $I$ in $V$,
such that ${\bf 1} \notin I$, $\omega \notin I$, the
quotient $V/I$ admits a natural vertex operator algebra structure.

Let $(M, Y_{M})$ be a weak module for a vertex operator algebra
$V$ (cf. [L]). A {\it ${\Z}_{+}$-graded weak} $V$-module ([FZ]) is a weak $V$-module
$M$ together with a ${\Z}_{+}$-gradation
$M=\oplus_{n=0}^{\infty}M(n)$ such that
\begin{eqnarray}
& &a_{m}M(n)\subseteq M(n+r-m-1)\;\;\;\mbox{for }a \in
V_{(r)},m,n,r\in {\Z},
\end{eqnarray}
where $M(n)=0$ for $n < 0$ by definition.

A weak $V$-module $M$ is called a {\it $V$-module} if $L(0)$
acts semisimply on $M$ with the decomposition into
$L(0)$-eigenspaces $M=\oplus_{\alpha\in {\C}}M_{(\alpha)}$ such
that for any $\alpha \in {\C}$, $\dim M_{(\alpha)}<\infty$
and $M_{(\alpha+n)}=0$ for $n\in {\Z}$ sufficiently small.

\subsection{Zhu's $A(V)$ theory}

Let $V$ be a vertex operator algebra. Following [Z], we
define bilinear maps $* :  V \times V
\to V$ and $\circ  :  V \times V \to V$ as follows.
For any homogeneous $a \in V$ and for any $b \in V$, let
\begin{eqnarray}
a \circ b=\Res_{z}\frac{(1+z)^{{\rm wt} a}}{z^{2}}Y(a,z)b \nonumber \\
a * b=\Res_{z}\frac{(1+z)^{{\rm wt} a}}{z}Y(a,z)b \nonumber
\end{eqnarray}
and extend to $V \times V \to V$ by linearity. Denote by
$O(V)$ the linear span of elements of the form $a \circ b$, and by
$A(V)$ the quotient space $V/O(V)$. For $a \in V$, denote by
$[a]$ the image of $a$ under the projection of $V$ onto
$A(V)$. The multiplication $*$ induces the multiplication
on $A(V)$ and $A(V)$ has a structure of an associative algebra.

\begin{prop}[{[FZ], Proposition 1.4.2}] Let $I$ be an ideal of $V$. Assume
${\bf 1} \notin I$, $\omega \notin I$. Then the associative
algebra $A(V/I)$ is isomorphic to $A(V)/A(I)$, where
$A(I)$ is the image of $I$ in $A(V)$.
\end{prop}

For any homogeneous $a \in V$ we define $o(a)=a_{{\rm wt} a-1}$
and extend this map linearly to $V$.

\begin{prop} [{[Z]}, Theorem 2.1.2, Theorem 2.2.1] \mbox{} \\
(a) Let $M=\oplus_{n=0}^{\infty}M(n)$ be a ${\Z}_{+}$-graded
weak $V$-module. Then $M(0)$ is an $A(V)$-module defined as
follows:
$$[a].v=o(a)v,$$
for any $a \in V$ and $v \in M(0)$. \\
(b) Let $U$ be an $A(V)$-module. Then there exists a ${\Z}_{+}$-graded
weak $V$-module $M$ such that the $A(V)$-modules $M(0)$ and $U$ are
isomorphic.
\end{prop}

\begin{prop}[{[Z], Theorem 2.2.2}] \label{t.1.3.5}
The equivalence classes of the irreducible $A(V)$-modules and
the equivalence classes of the irreducible ${\Z}_{+}$-graded weak $V$-modules
are in one-to-one correspondence.
\end{prop}

\subsection{Modules for affine Lie algebras}

Let ${\frak g}$ be a simple Lie algebra over ${\C}$ with a
triangular decomposition \linebreak ${\frak g}={\frak n}_{-} \oplus
{\frak h} \oplus {\frak n}_{+}$. Let $\Delta$ be the root system of
$({\frak g}, {\frak h})$, $\Delta_{+}\subset \Delta$
the set of positive roots, $\theta$ the highest root and
$(\cdot, \cdot): {\frak g}\times {\frak g}\to {\C}$ the Killing form,
normalized by the condition $(\theta, \theta)=2$.
Denote by $\Pi= \{ \alpha _{1}, \ldots , \alpha _{l} \}$ the set of simple
roots of $\frak g$, and by $\Pi ^{\vee}= \{ h_{1}, \ldots , h_{l} \}$
the set of simple coroots of $\frak g$.

The affine Lie algebra $\hat{\frak g}$ associated to ${\frak g}$
is the vector space ${\frak g}\otimes {\C}[t, t^{-1}] \oplus {\C}c
$ equipped with the usual bracket operation and the canonical central
element~$c$ (cf. [K]).
Let $h^{\vee}$ be the dual Coxeter number of $\hat{\frak g}$.
Let $\hat{\frak g}=\hat{\frak n}_{-} \oplus \hat{\frak h}
\oplus \hat{\frak n}_{+}$ be the corresponding
triangular decomposition of $\hat{\frak g}$.

Denote by $\hat{\Delta}$ the set of roots of $\hat{\frak g}$,
by $\hat{\Delta}_{+}$ the set of positive roots of $\hat{\frak g}$,
and by $\hat{\Pi}$ the set of simple
roots of $\hat{\frak g}$. Denote by $\hat{\Delta}^{\mbox{\scriptsize{re}}}$
the set of real roots of $\hat{\frak g}$ and let
$\hat{\Delta}^{\mbox{\scriptsize{re}}}_{+}= \hat{\Delta}^{\mbox{\scriptsize{re}}}
\cap \hat{\Delta}_{+}$. With $\alpha ^{\vee}$ denote the coroot of
a real root $\alpha \in \hat{\Delta}^{\mbox{\scriptsize{re}}}$.
For any $\lambda\in \hat{{\frak h}}^{*}$ set
$ D(\lambda)=$ $ \{ \lambda - \sum_{\alpha \in \hat{\Pi}}n_{\alpha}\alpha
\ \vert \ n_{\alpha} \in {\Z}_{+} \ \forall \alpha \in \hat{\Pi} \}.$

We say that a $\hat{\frak g}$-module $M$ is from the {\it category
$\mathcal{O}$} ([K]) if Cartan subalgebra $\hat{\frak h}$ acts semisimply on $M$ with
finite-dimensional weight spaces and there exists a finite
number of elements $\nu_{1}, \ldots , \nu_{k}
\in \hat{{\frak h}}^{*}$ such that for every weight
$\nu$ of $M$ holds $\nu \in \cup_{i=1}^{k}D(\nu_{i})$.

For every weight $\lambda\in \hat{{\frak h}}^{*}$, denote by
$M( \lambda)$ the Verma module for $\hat{\frak g}$ with highest
weight $\lambda$, and by $L( \lambda)$ the irreducible
$\hat{\frak g}$-module with highest
weight~$\lambda$.

Let $U$ be a ${\frak g}$-module, and let $k\in {\C}$.
Let $\hat{\frak g}_{+}={\frak g}\otimes t{\C}[t]$ act trivially on $U$ and
$c$ as the scalar multiplication operator $k$. Considering
$U$ as a ${\frak g}\oplus {\C}c
\oplus \hat{\frak g}_{+}$-module, we have
the induced $\hat{\frak g}$-module (so called {\it generalized Verma
module})
$$N(k,U)=U(\hat{\frak g})\otimes_{U({\frak g}\oplus {\C}c
\oplus \hat{\frak g}_{+})} U.$$

For a fixed $\mu \in {\frak h}^{*}$, denote by $V(\mu)$ the
irreducible highest-weight \linebreak ${\frak g}$-module with highest weight
$\mu$. Denote by $P_{+}$ the set of dominant integral weights of
$\frak g$, i.e. $P_{+}=\{ \, \mu \in {\frak h}^{*} \ \vert \
\mu (h_{i}) \in {\Z}_{+}, \mbox{ for } i=1, \ldots ,l   \}$.
Denote by $\omega_{1}, \ldots , \omega_{l} \in P_{+}$ the
fundamental weights of $\frak g$, defined by
$\omega_{i}(h_{j})= \delta _{ij}$ for all $i,j=1, \ldots ,l$.

We shall use the notation $N(k, \mu)$ to denote the $\hat{\frak
g}$-module $N(k,V(\mu))$. Denote by $J(k, \mu)$ the
maximal proper submodule of $N(k, \mu)$ and  $L(k, \mu)=N(k, \mu)
/J(k, \mu)$.

Denote by $\Lambda_{0} \in \hat{{\frak h}}^{*}$ the weight
defined by $\Lambda_{0}(c)=1$ and $\Lambda_{0}(h)=0$ for any
$h \in \frak h$. $N(k, \mu)$ is a highest-weight module
with highest weight $k \Lambda_{0}+ \mu$, and a quotient
of the Verma module $M(k \Lambda_{0}+ \mu)$.
$L(k, \mu)$ is the unique irreducible highest-weight
module with highest weight $k \Lambda_{0}+ \mu$, i.e.
$L(k, \mu) \cong L(k \Lambda_{0}+ \mu)$.

\subsection{Admissible modules for affine Lie algebras}

Let $\hat{\Delta}^{\vee \mbox{\scriptsize{re}}}$\
(resp $\hat{\Delta}^{\vee \mbox{\scriptsize{re}}}_{+}$)
$\subset \hat{\frak h}$ be the set of real
(resp positive real) coroots of $ \hat{\frak g}$. Fix $\lambda \in \hat{\frak h}^*$.
Let $\hat{\Delta}^{\vee \mbox{\scriptsize{re}}}_{\lambda}=
\{\alpha\in \hat{\Delta}^{\vee \mbox{\scriptsize{re}}} \ |\ \langle\lambda,\alpha\rangle\in{\Z} \}$,
$\hat{\Delta}^{\vee \mbox{\scriptsize{re}}}_{\lambda +}=\hat{\Delta}^{\vee \mbox{\scriptsize{re}}}_{\lambda}
\cap \hat{\Delta}^{\vee \mbox{\scriptsize{re}}}_{+}$,
$\hat{\Pi}^{\vee}$ the set of simple coroots in
$\hat{\Delta}^{\vee \mbox{\scriptsize{re}}}$ and
$\hat{\Pi}^{\vee}_{\lambda }=
\{\alpha \in \hat{\Delta}^{\vee \mbox{\scriptsize{re}}}_{\lambda +}\ \vert \ \alpha$
not equal to a sum of several coroots from $\hat{\Delta}^{\vee \mbox{\scriptsize{re}}}_{\lambda +} \}$.
Define $\rho$ in the usual way, and denote by $w. \lambda$
the "shifted" action of an element $w$ of the Weyl group of
$ \hat{\frak g}$.

Recall that a weight  $\lambda \in \hat{\frak h} ^* $
is called {\it admissible} (cf. [KW1], [KW2] and [W])
if the following properties are satisfied:
\begin{eqnarray*}
& &\langle\lambda + \rho,\alpha\rangle \notin -{\Z}_+ \mbox{ for all }
\alpha \in \hat{\Delta}^{\vee \mbox{\scriptsize{re}}}_{+}, \\
& &{\Q} \hat{\Delta}^{\vee \mbox{\scriptsize{re}}}_{\lambda}={\Q} \hat{\Pi}^{\vee} .
\end{eqnarray*}
The irreducible $\hat{\frak g}$-module $L(\lambda)$
is called admissible if the weight $\lambda \in \hat{\frak h} ^* $
is admissible.

We shall use the following results of V. Kac and M. Wakimoto:
\begin{prop}[{[KW1], Corollary 2.1}] \label{t.KW1}
Let $\lambda $ be an admissible weight. Then
$$L(\lambda)=\frac {M(\lambda)}{\sum_{\alpha \in \hat{\Pi}^{\vee}_{\lambda }}U (
\hat{\frak g}) v^{\alpha}}\ ,$$
where $v^{\alpha}\in M(\lambda)$ is a singular vector of weight $r_{\alpha}.
\lambda$, the highest weight vector of $M(r_{\alpha}.\lambda)=\ U(\hat{\frak g})
v^{\alpha}\subset M(\lambda)$.
\end{prop}

\begin{prop}[{[KW2], Theorem 4.1}] \label{t.KW2}
Let $M$ be a $\hat{\frak g}$-module from the category $\mathcal{O}$ such
that for any irreducible subquotient $L(\nu )$ the weight $\nu $ is
admissible. Then $M$ is completely reducible.
\end{prop}

\subsection{Vertex operator algebras $N(k,0)$ and $L(k,0)$, for
\\ $k \neq - h^{\vee}$}

Since $V(0)$ is the one-dimensional
trivial ${\frak g}$-module, it can be identified with~${\C}$.
Denote by ${\bf 1}=1 \otimes 1 \in N(k,0)$.
We note that $N(k,0)$ is spanned by the elements of the form
$x_{1}(-n_{1}-1)\cdots x_{m}(-n_{m}-1){\bf 1}$, where  $x_{1}, \dots, x_{m}\in
{\frak g}$ and $n_{1}, \dots, n_{m}\in {\Z}_{+}$, with $x(n)$ denoting
the representation image of
$x\otimes t^{n}$ for $x\in{\frak g}$ and $n\in {\Z}$.
Vertex operator map
$Y(\cdot, z): N(k,0) \to  (\mbox{\rm End}\;N(k,0))[[z, z^{-1}]]$
is uniquely determined by defining $Y({\bf 1}, z)$ to be the
identity operator on $N(k,0)$ and
$$Y(x(-1){\bf 1}, z)=\sum_{n\in {\Z}}x(n)z^{-n-1},$$
for $x\in {\frak g}$.
In the case that $k\ne - h^{\vee}$, $N(k,0)$ has a Virasoro element
\begin{equation}
\omega=\frac{1}{2(k+h^{\vee})}\sum_{i=1}^{\dim {\frak g}}
x^{i}(-1)^{2}{\bf 1},
\end{equation}
where $\{x^{i}\}_{i=1, \dots, \dim {\frak g}}$ is an arbitrary
orthonormal basis of ${\frak g}$ with respect to the form $(\cdot, \cdot)$.
In [FZ], the following results are proved:

\begin{prop}[{[FZ], Theorem 2.4.1}]
If $k\ne -h^{\vee}$, the quadruple \linebreak
$(N(k,0), Y, {\bf 1}, \omega)$ defined above is a vertex operator algebra.
\end{prop}

\begin{prop} [{[FZ], Theorem 3.1.1}]
The associative algebra $A(N(k,0))$
is canonically isomorphic to $U (\frak g ) $.
The isomorphism is given by $F:A(N(k,0)) \to U (\frak g )$
\begin{eqnarray*}
F([x_1(-n_1 -1)\cdots x_m(-n_m -1){\bf 1}])= (-1)^{n_1+\cdots +n_m}
x_m \cdots x_1,
\end{eqnarray*}
for any $x_1, \ldots ,x_m \in \frak g$ and any $n_1, \ldots ,n_m \in {\Z}_{+}$.
\end{prop}

Since every $\hat{\frak g}$-submodule of $N(k,0)$ is also an ideal
in the vertex operator algebra $N(k,0)$, it follows that
$L(k,0)$ is a vertex operator algebra, for any $k\ne -h^{\vee}$.
The associative algebra $A(L(k,0))$ is identified in
the next proposition, in the case when the maximal $\hat{\frak g}$-submodule of
$N(k,0)$ is generated by one singular vector.

\begin{prop} \label{p.1.5.5}
Assume that the maximal $\hat{\frak g}$-submodule of
$N(k,0)$ is generated by a singular vector, i.e.
$J(k,0)=U(\hat{\frak g})v_{sing}.$
Then
$$A(L(k,0)) \cong \frac{U(\frak g)}{\langle Q \rangle},$$
where $\langle Q \rangle$ is the two-sided ideal of $U(\frak g)$
generated by $Q=F([v_{sing}])$. \\ \noindent
Let $U$ be a $\frak g$-module. Then $U$ is an
$A(L(k,0))$-module if and only if $QU=0$.
\end{prop}

\section{Modules for vertex operator algebra \\
$L(n-l+\frac{1}{2},0)$
associated to affine Lie \\ algebra of type $B_{l}^{(1)}$}

Let $\frak g$ be the simple Lie algebra of type $B_{l}$,
and $\hat{\frak g}$ the affine Lie algebra associated to
$\frak g$. In this section we study the vertex operator algebra
\linebreak
$L(n-l+\frac{1}{2},0)$ associated to $\hat{\frak g}$,
for a positive integer $n$. Using results
from [FZ], [Z], we identify the corresponding associative algebra
$A(L(n-l+\frac{1}{2},0))$ with a quotient of $U({\frak g})$ by
the ideal generated by a certain vector. Using methods from
[MP], [A2], we show that the highest weights $\mu \in P_{+}$
of irreducible finite-dimensional $A(L(n-l+\frac{1}{2},0))$-modules
are characterized by the condition
$$( \mu , \epsilon_{1} ) \leq n-\frac{1}{2}, $$
where $\epsilon_{1}$ is the maximal short root for ${\frak g}$.
Using Zhu's theory we obtain the classification of irreducible
$L(n-l+\frac{1}{2},0)$-modules.
From this classification and results from [KW2], we obtain the
semisimplicity of the category of $L(n-l+\frac{1}{2},0)$-modules.
Using similar techniques, we also show that
there are finitely many irreducible weak
$L(n-l+\frac{1}{2},0)$-modules from the category $\mathcal{O}$.

\subsection{Simple Lie algebra of type $B_{l}$}

Let $\Delta =\{ \pm \epsilon_{i} \, \vert \ i=1, \ldots ,l \} \cup
\{ \pm (\epsilon_{i} \pm \epsilon_{j}) \, \vert \ i,j=1, \ldots ,l, \, i \neq j \}$
be the root system of type $B_{l}$. Fix the set of positive roots
$\Delta_{+}=\{ \epsilon_{i}, \, \vert \ i=1, \ldots ,l \} \cup$
$\{ \epsilon_{i} - \epsilon_{j} \, \vert \ i<j \} \cup
\{ \epsilon_{i} + \epsilon_{j} \, \vert \ i \neq j \}$.
Then the simple roots are $\alpha_{1}= \epsilon_{1} - \epsilon_{2}$,
$\alpha_{2}= \epsilon_{2} - \epsilon_{3}$, \ldots ,
$\alpha_{l-1}= \epsilon_{l-1} - \epsilon_{l}$,
$\alpha_{l}= \epsilon_{l}$.
The highest root is $\theta = \epsilon_{1} + \epsilon_{2}=
\alpha_{1}+2 \alpha_{2}+2 \alpha_{3}+ \cdots +2 \alpha_{l}$.

Let $\frak g$ be the simple Lie algebra associated to
the root system of type $B_{l}$. Let $e_{i},f_{i},h_{i}$,
$i=1, \ldots ,l$ be the Chevalley generators of $\frak g$.
Fix the root vectors:
\begin{eqnarray*}
e_{\epsilon_i - \epsilon_j} \!\!\!\!&=& \!\!\![e_{i},[e_{i+1},[ \ldots
[e_{j-2},e_{j-1}] \ldots] \, ] \, ] \qquad i<j \\
e_{\epsilon_i} \!\!\!\!&=&\!\!\! [e_{i},[e_{i+1},[ \ldots
[e_{l-1},e_{l}] \ldots] \, ] \, ]  \\
e_{\epsilon_i + \epsilon_j} \!\!\!\!&=&\!\!\! \frac{1}{2}[e_{\epsilon_i},
e_{\epsilon_j}] \qquad i<j \\
f_{\epsilon_i - \epsilon_j} \!\!\!\!&=& \!\!\! [f_{j-1},[f_{j-2},[ \ldots
[f_{i+1},f_{i}] \ldots] \, ] \, ] \qquad i<j \\
f_{\epsilon_i} \!\!\!\!&=& \!\!\! [f_{l},[f_{l-1},[ \ldots
[f_{i+1},f_{i}] \ldots] \, ] \, ]  \\
f_{\epsilon_i + \epsilon_j} \!\!\!\!&=& \!\!\! \frac{1}{2}[f_{\epsilon_j},
f_{\epsilon_i}] \qquad i<j
\end{eqnarray*}
Denote by $h_{\alpha}= \alpha ^{\vee}= [e_{\alpha},f_{\alpha}]$
coroots, for any positive root $\alpha \in \Delta_{+}$.
It is clear that $h_{\alpha_{i}}=h_{i}$. Let ${\frak g}={\frak n}_{-} \oplus
{\frak h} \oplus {\frak n}_{+}$ be the corresponding
triangular decomposition of ${\frak g}$.

\subsection{Vertex operator algebra $L(n-l+\frac{1}{2},0)$
associated to affine Lie algebra of type $B_{l}^{(1)}$}

Let $\hat{\frak g}$ be the affine Lie algebra associated to simple Lie
algebra of type $B_{l}$.
We want to show that the maximal $\hat{\frak g}$-submodule
of $N(n-l+\frac{1}{2},0)$ for $n \in {\N}$ is generated by a
singular vector. We need two lemmas to prove that.

Denote by $\lambda_{n}$ the weight $(n-l+\frac{1}{2})\Lambda_{0}$
for $n \in {\N}$.
Then $N(n-l+\frac{1}{2},0)$ is a quotient of $M(\lambda_{n})$
and $L(n-l+\frac{1}{2},0) \cong L(\lambda_{n})$.

\begin{lem} \label{l.2.3.1}
The weight $\lambda_{n}=(n-l+\frac{1}{2})\Lambda_{0}$
is admissible and
$$ \hat{\Pi}^{\vee}_{\lambda _{n}}= \{ (\delta - \epsilon_{1})^{\vee},
\alpha_{1}^{\vee},\alpha_{2}^{\vee}, \ldots ,
\alpha_{l}^{\vee} \},$$
for any $n \in {\N}$. Furthermore
\begin{eqnarray*}
&& \langle \lambda _{n}+ \rho,\alpha _{i}^{\vee}\rangle = 1 \ \ \mbox{for } i=1,
\ldots ,l, \\
&& \langle \lambda _{n} + \rho,(\delta - \epsilon_{1})^{\vee} \rangle =
2n.
\end{eqnarray*}
\end{lem}
{\bf Proof:} We have to show
\begin{eqnarray*}
& &\langle \lambda_{n} + \rho,\tilde{\alpha}^{\vee}\rangle \notin -{\Z}_{+} \mbox{ for
any }
\tilde{\alpha} \in \hat{\Delta}^{\mbox{\scriptsize{re}}}_{+}, \\
& &{\Q} \hat{\Delta}^{\vee \mbox{\scriptsize{re}}}_{\lambda _{n}}={\Q} \hat{\Pi}^{\vee} .
\end{eqnarray*}
Any positive real root $\tilde{\alpha} \in \hat{\Delta}^{\mbox{\scriptsize{re}}}_{+}$
of $\hat{\frak g}$
is of the form $\tilde{\alpha}= \alpha +m \delta$, for $m>0$ and $\alpha
\in \Delta$ or $m=0$ and $\alpha \in {\Delta}_{+}$.
Denote by $\bar{\rho}$ the sum of fundamental weights of
$\frak g$. Then $\rho =h^{\vee} \Lambda_{0}+ \bar{\rho}=(2l-1) \Lambda_{0}
+\bar{\rho}$, and
$\lambda_{n} + \rho = (n+l-\frac{1}{2})\Lambda_{0}+ \bar{\rho}.$
We have
\begin{eqnarray*}
\langle \lambda_{n} + \rho,\tilde{\alpha}^{\vee}\rangle &=&
\langle (n+l-\frac{1}{2})\Lambda_{0}+ \bar{\rho},(\alpha +m
\delta)^{\vee} \rangle \\
&=& \frac{2}{( \alpha , \alpha )}\left(m(n+l-\frac{1}{2})+( \bar{\rho}
, \alpha )\right).
\end{eqnarray*}
If $m=0$ and $\alpha \in {\Delta}_{+}$, then $( \bar{\rho}, \alpha
)>0$, and it is clear that
$\langle \lambda_{n} + \rho,\tilde{\alpha}^{\vee}\rangle \notin -{\Z}_{+}$.

Let $m \geq 1$. We have two cases: $( \alpha , \alpha
)=2$ and $( \alpha , \alpha )=1$.

If $( \alpha , \alpha )=2$ and $m$ is odd, then
$\langle \lambda_{n} + \rho,\tilde{\alpha}^{\vee}\rangle \notin {\Z}$, which implies
\linebreak
$\langle \lambda_{n} + \rho,\tilde{\alpha}^{\vee}\rangle \notin -{\Z}_{+}$.

If $( \alpha , \alpha )=2$ and $m$ is even, then $m \geq 2$,
and since $( \bar{\rho}, \alpha ) \geq -(2l-1)$ for any
$\alpha \in \Delta$, we have
\begin{eqnarray*}
\langle \lambda_{n} + \rho,\tilde{\alpha}^{\vee}\rangle=
m(n+l-\frac{1}{2})+( \bar{\rho}, \alpha )\geq
2(n+l-\frac{1}{2})-(2l-1)=2n \geq 2,
\end{eqnarray*}
which implies
$\langle \lambda_{n} + \rho,\tilde{\alpha}^{\vee}\rangle \notin -{\Z}_{+}$.

If $( \alpha , \alpha )=1$, then $\alpha = \pm
\epsilon_{i}$ for some $i=1, \ldots ,l$, which implies
$( \bar{\rho}, \alpha ) \geq -l$. Then
\begin{eqnarray*}
\langle \lambda_{n} + \rho,\tilde{\alpha}^{\vee}\rangle=
2 \left(m(n+l-\frac{1}{2})+( \bar{\rho} , \alpha )\right) \geq
2 (n+l-\frac{1}{2}-l)=2n-1 \geq 1,
\end{eqnarray*}
which implies
$\langle \lambda_{n} + \rho,\tilde{\alpha}^{\vee}\rangle \notin -{\Z}_{+}$.

Thus,
$\langle \lambda_{n} + \rho,\tilde{\alpha}^{\vee}\rangle \notin -{\Z}_{+} \mbox{ for
any } \tilde{\alpha} \in \hat{\Delta}^{\mbox{\scriptsize{re}}}_{+}$.
Furthermore,
\begin{eqnarray*}
&& \langle \lambda _{n}+ \rho,\alpha _{i}^{\vee}\rangle = 1 \ \ \mbox{for } i=1,
\ldots ,l, \\
&& \langle \lambda _{n} + \rho,(\delta - \epsilon_{1})^{\vee} \rangle =
2n \in \N,
\end{eqnarray*}
which implies
$$ \hat{\Pi}^{\vee}_{\lambda _{n}}= \{ (\delta - \epsilon_{1})^{\vee},
\alpha_{1}^{\vee},\alpha_{2}^{\vee}, \ldots ,
\alpha_{l}^{\vee} \}.$$
It follows that
${\Q} \hat{\Delta}^{\vee \mbox{\scriptsize{re}}}_{\lambda _{n} }=
{\Q} \hat{\Pi}^{\vee}_{\lambda _{n}}=
{\Q} \hat{\Pi}^{\vee}$. $\;\;\;\;\Box$

\begin{lem} \label{l.2.3.2}
Vector
$$v_n=(-\frac{1}{4}e_{\epsilon_1}(-1)^2+ e_{\epsilon_1 - \epsilon_2}(-1)
e_{\epsilon_1 + \epsilon_2}(-1)+ \cdots +e_{\epsilon_1 - \epsilon_l}(-1)
e_{\epsilon_1 + \epsilon_l}(-1))^n{\bf 1} $$
is a singular vector in $N(n-l+\frac{1}{2},0)$ for any $n \in {\N}$.
\end{lem}
{\bf Proof:} It is sufficient to show
\begin{eqnarray*}
& & e_{i}(0).v_n=0, \ i=1, \ldots ,l \\
& & f_{\theta}(1).v_n=0.
\end{eqnarray*}
We introduce the notation
$$ u=-\frac{1}{4}e_{\epsilon_1}(-1)^2+ e_{\epsilon_1 - \epsilon_2}(-1)
e_{\epsilon_1 + \epsilon_2}(-1)+ \cdots +e_{\epsilon_1 - \epsilon_l}(-1)
e_{\epsilon_1 + \epsilon_l}(-1) \in U(\hat{\frak g}). $$
It can easily be checked that vectors $e_{i}(0)$ commute with $u$ in $U(\hat{\frak g})$
for any $i=1, \ldots ,l$, which implies
\begin{eqnarray*}
e_{i}(0).v_n=e_{i}(0).u^{n}{\bf 1}=[e_{i}(0),u^{n}]{\bf 1}=0, \quad i=1, \ldots
,l.
\end{eqnarray*}
Similarly, we can show
\begin{eqnarray*}
[f_{\theta}(1),u]&=&(k+l-\frac{3}{2})e_{1}(-1)-
e_{1}(-1)h_{\epsilon_1 + \epsilon_2}(0)-
\frac{1}{2} e_{\epsilon_1}(-1)f_{\epsilon_2}(0) \\
& - & e_{\epsilon_1 + \epsilon_3}(-1)f_{\epsilon_2 +
\epsilon_3}(0) -
e_{\epsilon_1 - \epsilon_3}(-1)f_{\epsilon_2 - \epsilon_3}(0)-
\ldots \\
& - & e_{\epsilon_1 + \epsilon_l}(-1)f_{\epsilon_2 + \epsilon_l}(0)
-e_{\epsilon_1 - \epsilon_l}(-1)f_{\epsilon_2 - \epsilon_l}(0).
\end{eqnarray*}
It can easily be checked that vectors
$f_{\epsilon_2}(0),f_{\epsilon_2 + \epsilon_3}(0),f_{\epsilon_2 - \epsilon_3}(0),
\ldots ,f_{\epsilon_2 + \epsilon_l}(0),f_{\epsilon_2 -
\epsilon_l}(0)$ commute with $u$ in $U(\hat{\frak g})$ and that
\begin{eqnarray*}
[h_{\epsilon_1 + \epsilon_2}(0),u^{m}]=2mu^{m}, \ \mbox{for } m \in {\N}.
\end{eqnarray*}
We have:
\begin{eqnarray*}
f_{\theta}(1).v_n &=&f_{\theta}(1).u^{n}{\bf 1}=[f_{\theta}(1),u^{n}]{\bf
1}=[f_{\theta}(1),u]u^{n-1}{\bf 1}+u[f_{\theta}(1),u]u^{n-2}{\bf
1}+ \cdots \\
&+& u^{n-1}[f_{\theta}(1),u]{\bf 1}=(k+l-\frac{3}{2})e_{1}(-1)u^{n-1}{\bf 1}
-2(n-1)e_{1}(-1)u^{n-1}{\bf 1} \\
&+& (u(k+l-\frac{3}{2})e_{1}(-1)u^{n-2}
-ue_{1}(-1) \cdot 2(n-2)u^{n-2}){\bf 1}+ \cdots \\
&+& u^{n-1}(k+l-\frac{3}{2})e_{1}(-1){\bf 1} \\
&= &n(k+l-\frac{3}{2})e_{1}(-1)u^{n-1}{\bf 1}-2((n-1)+(n-2)+ \cdots
+1)e_{1}(-1)u^{n-1}{\bf 1} \\
&= &(n(k+l-\frac{3}{2})-2 \cdot \frac{(n-1)n}{2})e_{1}(-1)u^{n-1}{\bf
1} \\
&= &n(k+l-n-\frac{1}{2})e_{1}(-1)u^{n-1}{\bf 1}=0,
\end{eqnarray*}
because $k=n-l+\frac{1}{2}$.
Thus $\hat{\frak n}_{+}.v_n=0$, and $v_n$ is a singular vector in
$N(n-l+\frac{1}{2},0)$.  $\;\;\;\;\Box$

\begin{thm} \label{t.2.3.3} The maximal $\hat{\frak g}$-submodule
of $N(n-l+\frac{1}{2},0)$ is $J(n-l+\frac{1}{2},0)=
U(\hat{\frak g})v_{n}$, where
$$v_n=(-\frac{1}{4}e_{\epsilon_1}(-1)^2+ e_{\epsilon_1 - \epsilon_2}(-1)
e_{\epsilon_1 + \epsilon_2}(-1)+ \cdots +e_{\epsilon_1 - \epsilon_l}(-1)
e_{\epsilon_1 + \epsilon_l}(-1))^n{\bf 1}, \ n \in {\N}. $$
\end{thm}
{\bf Proof:} It follows from Proposition \ref{t.KW1} and Lemma \ref{l.2.3.1}
that the maximal submodule of the Verma module $M(\lambda_{n})$ is generated by $l+1$
singular vectors with weights
$$r_{\delta - \epsilon_1}.\lambda_{n}, r_{\alpha_{1}}.\lambda_{n},
\ldots ,r_{\alpha_{l}}.\lambda_{n}.$$
It follows from Lemma \ref{l.2.3.2} that $v_n$ is a singular vector of
weight
$\lambda_{n}-2n \delta +2n \epsilon_1=r_{\delta -
\epsilon_1}.\lambda_{n}$. Other singular vectors
have weights
$$r_{\alpha_{i}}.\lambda_{n}=\lambda_{n}- \langle \lambda_{n}+
\rho, \alpha_{i}^{\vee} \rangle \alpha_{i}=\lambda_{n}
- \alpha_{i}, \ i=1,\ldots ,l,$$
so the images of these vectors under the projection of
$M(\lambda_{n})$ onto \linebreak $N(n-l+\frac{1}{2},0)$ are 0.
Therefore, the maximal submodule of $N(n-l+\frac{1}{2},0)$ is
generated by the vector $v_n$, i.e. $J(n-l+\frac{1}{2},0)=
U(\hat{\frak g})v_{n}$. $\;\;\;\;\Box$

Using Theorem \ref{t.2.3.3} and Proposition \ref{p.1.5.5}
we can identify the associative algebra $A(L(n-l+\frac{1}{2},0))$:

\begin{prop} \label{t.2.3.5}
The associative algebra $A(L(n-l+\frac{1}{2},0))$ is
isomorphic to the algebra $U(\frak g)/I_{n}$, where $I_{n}$
is the two-sided ideal of $U(\frak g)$ generated by
$$v_n'=(-\frac{1}{4}e_{\epsilon_1}^2+ e_{\epsilon_1 - \epsilon_2}
e_{\epsilon_1 + \epsilon_2}+ \cdots +e_{\epsilon_1 - \epsilon_l}
e_{\epsilon_1 + \epsilon_l})^n.$$
\end{prop}

\subsection{Modules for associative algebra $A(L(n-l+\frac{1}{2},0))$}
\label{subsec.3.3}

In this subsection we present the method from [MP], [A2], [AM] for
classification of irreducible $A(L(n-l+\frac{1}{2},0))$-modules
from the category $\mathcal{O}$ by solving certain systems of
polynomial equations.

Denote by $_L$ the adjoint action of  $U(\frak g)$ on
$U(\frak g)$ defined by $ X_Lf=[X,f]$ for $X \in
\frak g$ and $f \in U(\frak g)$. Let $R$ be a $U(\frak g)$-submodule
of $U(\frak g)$ generated by the vector $v_n'$.
Clearly, $R$ is an irreducible finite-dimensional $U(\frak g)$-module
isomorphic to $V(2n \epsilon_{1})$. Let $R_{0}$ be the zero-weight
subspace of $R$.

\begin{prop}[{[A2], Proposition 2.4.1}, {[AM], Lemma 3.4.3}]
Let $V(\mu)$ be an irreducible highest weight $U(\frak g)$-module
with the highest weight vector $v_{\mu}$, for $\mu \in {\frak h}^{*}$.
The following statements are equivalent: \\
$\, \!$ (1) $V(\mu)$ is an $A(L(n-l+\frac{1}{2},0))$-module, \\
$\, \!$ (2) $RV(\mu)=0$, \\
$\, \!$ (3) $R_{0}v_{\mu}=0.$
\end{prop}
{\bf Proof:} The equivalence of (1) and (2) follows from the fact that
$R$ generates the ideal $I_{n}$. Clearly, (2) implies (3). To prove
the converse, suppose that $R_{0}v_{\mu}=0.$
We claim that
$RV(\mu)$ is a $\frak g$-submodule of $V(\mu)$.
Let $x \in \frak g$, $r \in R$ and $v \in V(\mu)$. We get
\[ x(rv)=[x,r]v + r(xv). \]
Since $R$ is a $\frak g$-module, $[x,r] \in R$, and that
implies $x(rv) \in RV(\mu)$. Since $V(\mu)$ is irreducible,
$RV(\mu)=0$ or $RV(\mu)=V(\mu)$. To prove that $RV(\mu)=0$,
it is enough to show that $v_{\mu} \notin RV(\mu)$.
Clearly,
\[ RV(\mu)= R U({\frak n}_{-})v_{\mu}= U({\frak n}_{-})Rv_{\mu},\]
since $R$ is a $\frak g$-module under the adjoint action.
Since $R \subseteq U(\frak g)$, the Poincar\'{e}-Birkhoff-Witt
theorem implies that every
element of $r \in R$ can be written as a linear combination
of elements of the form $r_{0}r_{-}r_{+}$, where
$r_{0} \in S( \frak h)$, $r_{-} \in U({\frak n}_{-})$ and
$r_{+} \in U({\frak n}_{+})$. It follows that, if the weight of
$r$ is positive, then $rv_{\mu}=0$, and if the weight of $r$
is negative, then the weight of $rv_{\mu}$ is $\mu + wt(r) <\mu$.
From this we obtain that $v_{\mu} \in RV(\mu)$ if and only if
$R_{0}v_{\mu}\neq 0.$ Thus, $R_{0}v_{\mu}=0$ implies
$v_{\mu} \notin RV(\mu)$, which implies $RV(\mu)=0$. $\;\;\;\;\Box$

Let $r \in R_{0}$. Clearly there exists the unique polynomial
$p_{r} \in S( \frak h)$ such that
$$ rv_{\mu}=p_{r}(\mu)v_{\mu}.$$
Set $ {\mathcal P}_{0}=\{ \ p_{r} \ \vert \ r \in R_{0} \}.$
We have:

\begin{coro} \label{c.1.7.2} There is one-to-one correspondence between \\
$\, \!$ (1) irreducible $A(L(n-l+\frac{1}{2},0))$-modules
from the category $\mathcal{O}$, \\
$\, \!$ (2) weights $\mu \in {\frak h}^{*}$ such that
$p(\mu)=0$ for all $p \in {\mathcal P}_{0}$.
\end{coro}

\subsection{Construction of some polynomials in $ {\mathcal P}_{0}$}

The following lemmas are obtained by direct
calculations in $U(\frak g)$:

\begin{lem} \label{l.2.4.0}
Let $X \in \frak g$ and $Y_{1}, \ldots ,
Y_{m} \in U(\frak g)$. Then
$$(X^{n})_{L}(Y_{1} \cdot \ldots \cdot Y_{m})=
{\displaystyle \sum_{{(k_{1},\ldots ,k_{m}) \in {\Z}_{+}^{m}
\atop \sum k_{i}=n}}} {n \choose k_{1},\ldots ,k_{m}}
(X^{k_{1}})_{L}Y_{1} \cdot \ldots \cdot (X^{k_{m}})_{L}Y_{m}.$$
\end{lem}

\begin{lem} \label{l.2.4.1}
\begin{eqnarray*}
& &(1) \ (e_{\alpha}^{m}) _L (f_{\alpha}^{m}) \in
m!\cdot h_{\alpha} \cdot \ldots \cdot (h_{\alpha} -m+1)
+U(\frak g) e_{\alpha}, \ \forall \alpha \in \Delta_{+};  \\
& &(2)\ (e_{\alpha}^{k}) _L (f_{\alpha}^{m}) \in
U(\frak g)e_{\alpha}, \ \mbox{ for } k>m \mbox{ and }
\ \forall \alpha \in \Delta_{+};   \\
& &(3)\ (e_{\epsilon_1}^{2k})_L (f_{\epsilon_1 + \epsilon_i }^{k})
=(-1)^{k}(2k)! \cdot e_{\epsilon_1 - \epsilon_i }^k, \
 i=2,\ldots,l;  \\
& &(4)\ (e_{\epsilon_1}^{2k+j})_L (f_{\epsilon_1 + \epsilon_i }^k)
=0, \ \mbox{ for } j>0 \mbox{ and } \  i=2,\ldots,l;  \\
& &(5)\ (e_{\epsilon_1}^{r})_L (f_{\epsilon_1 - \epsilon_i }^k)
\in U(\frak g){\frak n}_{+}, \ \mbox{ for } r>0 \mbox{ and } \  i=2,\ldots,l;  \\
& &(6)\ e_{\alpha}^kp(h)= p(h-k\alpha(h))e_{\alpha}^k, \ \forall p \in S({\frak h}); \\
& &(7)\ (e_{\epsilon_1 +\epsilon_{i}}^{k})_L (f_{\epsilon_1 + \epsilon_2 }^{m})
=m(m-1) \cdots (m-k+1) f_{\epsilon_1 + \epsilon_2 }^{m-k}
f_{\epsilon_2 - \epsilon_{i} }^{k},  \\
&& \qquad  \mbox{for }i=3,\ldots,l \ \mbox{and } k \leq m; \\
& &(8)\ (e_{\epsilon_1 +\epsilon_{i}}^{k})_L (f_{\epsilon_1 - \epsilon_2 }^{m})
\in U(\frak g){\frak n}_{+},  \  \mbox{for }i=3,\ldots,l \ \mbox{and } k >0; \\
& &(9)\ (e_{\epsilon_1 +\epsilon_{2}}^{k})_L (f_{\epsilon_2 - \epsilon_{i} }^{m})
\in U(\frak g)e_{\epsilon_1 +\epsilon_{i}},  \  \mbox{for }i=3,\ldots,l \ \mbox{and } k >0; \\
&&(10) \ (e_{\alpha}^{k}) _L (f_{\alpha}^{m}) \in
m(m-1) \cdots (m-k+1)f_{\alpha}^{m-k} \cdot \\
&& \qquad \ \cdot (h_{\alpha}-m+k) \cdot \ldots \cdot (h_{\alpha} -m+1)
+U(\frak g) e_{\alpha}, \ \forall \, \alpha \in \Delta_{+} \ \mbox{and } k \leq m; \\
& &(11)\ (e_{\epsilon_1 -\epsilon_{i}}^{k})_L (f_{\epsilon_2 - \epsilon_{i} }^{k})
=k! \, e_{\epsilon_1 - \epsilon_2 }^{k} \ \mbox{for } i=3,\ldots,l; \\
& &(12)\ (e_{\epsilon_1 -\epsilon_{i}}^{k})_L (f_{\epsilon_1 - \epsilon_2 }^{m})
\in U(\frak g){\frak n}_{+},  \  \mbox{for }i=3,\ldots,l \ \mbox{and } k>0.
\end{eqnarray*}
\end{lem}

\begin{lem}\label{l.2.4.2.1}
Let $\beta _{1}, \ldots ,\beta _{k}, \gamma _{1}, \ldots ,\gamma _{m}
\in \Delta_{+}$ such that $\sum _{i=1}^{k} \beta _{i}= \sum _{i=1}^{m}
\gamma _{i}$. Let $Y_{1},Y_{2} \in U(\frak g)$ such that
\begin{eqnarray*}
&&Y_{1}=e_{\beta _{1}} \cdots e_{\beta _{k}}, \quad
[e_{\beta _{i}},e_{\beta _{j}}]=0, \ \mbox{for all }\, i,j, \\
&&Y_{2}=f_{\gamma _{1}} \cdots f_{\gamma _{m}}, \quad
[f_{\gamma _{i}},f_{\gamma _{j}}]=0, \ \mbox{for all }\, i,j.
\end{eqnarray*}
Then
\begin{eqnarray*}
&&Y_{1} \,{}_{L}Y_{2} \in e_{\beta _{1}} \cdots e_{\beta _{k}}
f_{\gamma _{1}} \cdots f_{\gamma _{m}}+  U(\frak g){\frak n}_{+}, \\
&&Y_{2} \,{}_{L}Y_{1} \in (-1)^{m} e_{\beta _{1}} \cdots e_{\beta _{k}}
f_{\gamma _{1}} \cdots f_{\gamma _{m}}+  U(\frak g){\frak n}_{+}.
\end{eqnarray*}
\end{lem}

In the following lemma we calculate $l+1$ polynomials in the set
${\mathcal P}_{0}$. This is a crucial technical result needed for
classification of irreducible $L(n-l+\frac{1}{2},0)$-modules.

\begin{lem}\label{l.2.4.2}
Let
\begin{eqnarray*}
&\!\!\!\!\!\! (1)&\!\! q(h)= \!\!\!\!\!\!\! {\displaystyle
\sum_{{(k_{1},\ldots ,k_{l}) \in {\Z}_{+}^{l} \atop \sum k_{i}=n}}} \frac{1}{k_{1}!4^{k_{1}}} \cdot
(h_{\epsilon_1}-2k_{2}-\ldots -2k_{l})\cdot \ldots \cdot
(h_{\epsilon_1}-2n+1)\cdot \\
& & \ \ \ \ \ \ \ \ \ \cdot (h_{\epsilon_1 - \epsilon_l}
-k_{l-1}-\ldots -k_{2}) \cdot \ldots \cdot (h_{\epsilon_1 - \epsilon_l}
-k_{l-1}-\ldots -k_{2}-k_{l}+1) \cdot \\
& & \ \ \ \ \ \ \ \ \ \cdot \ldots \cdot
h_{\epsilon_1 - \epsilon_2} \cdot \ldots \cdot (h_{\epsilon_1 -
\epsilon_2}-k_{2}+1), \\
&\!\!\!\!\!\! (2)&\!\! p_{i}(h)=h_{i}(h_{i}-1)\cdot \ldots \cdot (h_{i}-n+1)
(h_{\epsilon_i + \epsilon_{i+1}}+l-i-\frac{1}{2})\cdot \\
&\!\!\!\!\!\!& \cdot (h_{\epsilon_i + \epsilon_{i+1}}+l-i-\frac{3}{2})\cdot \ldots
\cdot (h_{\epsilon_i + \epsilon_{i+1}}+l-n-i+\frac{1}{2}), \
\mbox{for } i=1, \ldots ,l-1, \\
&\!\!\!\!\!\! (3)&\!\! p_{l}(h)=h_{l}(h_{l}-1)\cdot \ldots \cdot (h_{l}-2n+1).
\end{eqnarray*}
Then $p_{1}, \ldots ,p_{l},q \in {\mathcal P}_{0}$.
\end{lem}
{\bf Proof:} (1) We claim that
\begin{eqnarray*}
(e_{\epsilon_1}^{2n}f_{\epsilon_1}^{4n})_L
(-\frac{1}{4}e_{\epsilon_1}^2+ e_{\epsilon_1 - \epsilon_2}
e_{\epsilon_1 + \epsilon_2}+ \cdots +e_{\epsilon_1 - \epsilon_l}
e_{\epsilon_1 + \epsilon_l})^n \in  c \, q(h)  +  U(\frak g){\frak n}_{+},
\end{eqnarray*}
for some $c \neq 0$.

We introduce the notation:
$$ \bar{u}= -\frac{1}{4}e_{\epsilon_1}^2+ e_{\epsilon_1 - \epsilon_2}
e_{\epsilon_1 + \epsilon_2}+ \cdots +e_{\epsilon_1 - \epsilon_l}
e_{\epsilon_1 + \epsilon_l} \in U(\frak g). $$
It can easily be checked that
$$(f_{\epsilon_1}^{4})_{L} \bar{u}=
24(-\frac{1}{4}f_{\epsilon_1}^2+ f_{\epsilon_1 - \epsilon_2}
f_{\epsilon_1 + \epsilon_2}+ \cdots +f_{\epsilon_1 - \epsilon_l}
f_{\epsilon_1 + \epsilon_l})$$
and $(f_{\epsilon_1}^{5})_{L} \bar{u}=0$,
which implies
\begin{eqnarray*}
& &(f_{\epsilon_1}^{4n})_L
(-\frac{1}{4}e_{\epsilon_1}^2+ e_{\epsilon_1 - \epsilon_2}
e_{\epsilon_1 + \epsilon_2}+ \cdots +e_{\epsilon_1 - \epsilon_l}
e_{\epsilon_1 + \epsilon_l})^n=
(f_{\epsilon_1}^{4n})_L(\bar{u}^{n}) \\
& &= \frac{(4n)!}{(4!)^{n}}
((f_{\epsilon_1}^{4})_L\bar{u})^{n}
 =(4n)! (-\frac{1}{4}f_{\epsilon_1}^2+ f_{\epsilon_1 - \epsilon_2}
f_{\epsilon_1 + \epsilon_2}+ \cdots +f_{\epsilon_1 - \epsilon_l}
f_{\epsilon_1 + \epsilon_l})^{n}.
\end{eqnarray*}
Therefore,
\begin{eqnarray} \label{2.4.2}
& &(e_{\epsilon_1}^{2n}f_{\epsilon_1}^{4n})_L
(-\frac{1}{4}e_{\epsilon_1}^2+ e_{\epsilon_1 - \epsilon_2}
e_{\epsilon_1 + \epsilon_2}+ \cdots +e_{\epsilon_1 - \epsilon_l}
e_{\epsilon_1 + \epsilon_l})^n \nonumber \\
& & = c'(e_{\epsilon_1}^{2n})_L
(-\frac{1}{4}f_{\epsilon_1}^2+ f_{\epsilon_1 - \epsilon_2}
f_{\epsilon_1 + \epsilon_2}+ \cdots +f_{\epsilon_1 - \epsilon_l}
f_{\epsilon_1 + \epsilon_l})^{n},
\end{eqnarray}
where $c' \neq 0$.
Since all root vectors $f_{\epsilon_1},
f_{\epsilon_1 - \epsilon_2}, f_{\epsilon_1 + \epsilon_2}, \ldots ,
f_{\epsilon_1 - \epsilon_l}$, $f_{\epsilon_1 + \epsilon_l} \in \frak g$
commute, we have
\begin{eqnarray} \label{2.4.1}
& &(-\frac{1}{4}f_{\epsilon_1}^2+ f_{\epsilon_1 - \epsilon_2}
f_{\epsilon_1 + \epsilon_2}+ \cdots +f_{\epsilon_1 - \epsilon_l}
f_{\epsilon_1 + \epsilon_l})^{n}  \\
& &= {\displaystyle \sum_{{(k_{1},\ldots ,k_{l}) \in {\Z}_{+}^{l} \atop \sum k_{i}=n}}} {n \choose k_{1},\ldots ,k_{l}}
(-1)^{k_{1}} \frac{1}{4^{k_{1}}} f_{\epsilon_1}^{2k_{1}}
f_{\epsilon_1 - \epsilon_2}^{k_{2}} f_{\epsilon_1 + \epsilon_2}^{k_{2}}
\cdots f_{\epsilon_1 - \epsilon_l}^{k_{l}} f_{\epsilon_1 +
\epsilon_l}^{k_{l}}. \nonumber
\end{eqnarray}
We calculate the action of $e_{\epsilon_1}^{2n}$ on every
summand above. By using claim $(5)$ from Lemma \ref{l.2.4.1}
we obtain
$$ (e_{\epsilon_1}^{r})_L (f_{\epsilon_1 - \epsilon_2 }^{k_{2}}
\cdots  f_{\epsilon_1 - \epsilon_l }^{k_{l}})
\in U(\frak g){\frak n}_{+}, \quad r>0,$$
which implies
\begin{eqnarray} \label{2.4.0}
& &(e_{\epsilon_1}^{2n})_L (f_{\epsilon_1}^{2k_{1}}
f_{\epsilon_1 - \epsilon_2}^{k_{2}} f_{\epsilon_1 + \epsilon_2}^{k_{2}}
\cdots f_{\epsilon_1 - \epsilon_l}^{k_{l}} f_{\epsilon_1 +
\epsilon_l}^{k_{l}}) \nonumber \\
& =&  (e_{\epsilon_1}^{2n})_L (f_{\epsilon_1}^{2k_{1}}
f_{\epsilon_1 + \epsilon_2}^{k_{2}} \cdots f_{\epsilon_1 +
\epsilon_l}^{k_{l}} f_{\epsilon_1 - \epsilon_2}^{k_{2}}
\cdots f_{\epsilon_1 - \epsilon_l}^{k_{l}})  \nonumber \\
& \in & [(e_{\epsilon_1}^{2n})_L (f_{\epsilon_1}^{2k_{1}}
f_{\epsilon_1 + \epsilon_2}^{k_{2}} \cdots f_{\epsilon_1 +
\epsilon_l}^{k_{l}})] (f_{\epsilon_1 - \epsilon_2}^{k_{2}}
\cdots f_{\epsilon_1 - \epsilon_l}^{k_{l}})+ U(\frak g){\frak
n}_{+}.
\end{eqnarray}
Using Lemma \ref{l.2.4.0} we can calculate:
\begin{eqnarray*}
& & (e_{\epsilon_1}^{2n})_L (f_{\epsilon_1}^{2k_{1}}
f_{\epsilon_1 + \epsilon_2}^{k_{2}} \cdots f_{\epsilon_1 +
\epsilon_l}^{k_{l}})=
(e_{\epsilon_1}^{2n})_L (f_{\epsilon_1 + \epsilon_2}^{k_{2}}
\cdots f_{\epsilon_1 + \epsilon_l}^{k_{l}}f_{\epsilon_1}^{2k_{1}})
\\
& = & {\displaystyle \sum_{{(m_{1},\ldots ,m_{l}) \in {\Z}_{+}^{l} \atop \sum m_{i}=2n}}} {2n \choose m_{1},\ldots ,m_{l}}
[(e_{\epsilon_1}^{m_{2}})_L (f_{\epsilon_1 + \epsilon_2}^{k_{2}})
\cdots (e_{\epsilon_1}^{m_{l}})_L (f_{\epsilon_1 + \epsilon_l}^{k_{l}})
(e_{\epsilon_1}^{m_{1}})_L(f_{\epsilon_1}^{2k_{1}})].
\end{eqnarray*}
It follows from claims $(2)$ and $(4)$ from Lemma \ref{l.2.4.1}
that only nontrivial summand above is for
$m_{1}=2k_{1}, m_{2}=2k_{2}, \ldots , m_{l}=2k_{l}$.
By using claims $(1)$ and $(3)$ from Lemma \ref{l.2.4.1}
we obtain:
\begin{eqnarray*}
& & (e_{\epsilon_1}^{2n})_L (f_{\epsilon_1}^{2k_{1}}
f_{\epsilon_1 + \epsilon_2}^{k_{2}} \cdots f_{\epsilon_1 +
\epsilon_l}^{k_{l}}) \\
& \in &\frac{(2n)!}{(2k_{1})! \cdots (2k_{l})!}
(e_{\epsilon_1}^{2k_{2}})_L (f_{\epsilon_1 + \epsilon_2}^{k_{2}})
\cdots (e_{\epsilon_1}^{2k_{l}})_L (f_{\epsilon_1 + \epsilon_l}^{k_{l}})
(e_{\epsilon_1}^{2k_{1}})_L(f_{\epsilon_1}^{2k_{1}}) + U(\frak g) e_{\epsilon_1} \\
& \in &(2n)!(-1)^{n-k_{1}} e_{\epsilon_1 - \epsilon_2}^{k_{2}}
\cdots e_{\epsilon_1 - \epsilon_l}^{k_{l}} \cdot h_{\epsilon_1}
\cdots (h_{\epsilon_{1}} -2k_{1}+1) + U(\frak g) e_{\epsilon_1}.
\end{eqnarray*}
By putting the expression above in relation (\ref{2.4.0}), we obtain:
\begin{eqnarray*}
& &(e_{\epsilon_1}^{2n})_L (f_{\epsilon_1}^{2k_{1}}
f_{\epsilon_1 - \epsilon_2}^{k_{2}} f_{\epsilon_1 + \epsilon_2}^{k_{2}}
\cdots f_{\epsilon_1 - \epsilon_l}^{k_{l}} f_{\epsilon_1 +
\epsilon_l}^{k_{l}}) \\
& \in & (2n)!(-1)^{n-k_{1}} e_{\epsilon_1 - \epsilon_2}^{k_{2}}
\cdots e_{\epsilon_1 - \epsilon_l}^{k_{l}} \cdot h_{\epsilon_1}
\cdots (h_{\epsilon_{1}} -2k_{1}+1)
f_{\epsilon_1 - \epsilon_2}^{k_{2}}
\cdots f_{\epsilon_1 - \epsilon_l}^{k_{l}}+ U(\frak g){\frak
n}_{+}
\end{eqnarray*}
By using claims $(6)$ and $(1)$ from Lemma \ref{l.2.4.1},
we have:
\begin{eqnarray*}
& &(e_{\epsilon_1}^{2n})_L (f_{\epsilon_1}^{2k_{1}}
f_{\epsilon_1 - \epsilon_2}^{k_{2}} f_{\epsilon_1 + \epsilon_2}^{k_{2}}
\cdots f_{\epsilon_1 - \epsilon_l}^{k_{l}} f_{\epsilon_1 +
\epsilon_l}^{k_{l}}) \\
& \in & (2n)!(-1)^{n-k_{1}}(h_{\epsilon_{1}}-2k_{2}-\ldots
-2k_{l}) \cdots (h_{\epsilon_{1}} -2n+1) \\
& & \cdot e_{\epsilon_1 - \epsilon_2}^{k_{2}}
\cdots e_{\epsilon_1 - \epsilon_l}^{k_{l}}
f_{\epsilon_1 - \epsilon_l}^{k_{l}}
\cdots f_{\epsilon_1 - \epsilon_2}^{k_{2}}+ U(\frak g){\frak
n}_{+} \\
& \in & (2n)!(-1)^{n-k_{1}}(h_{\epsilon_{1}}-2k_{2}-\ldots
-2k_{l}) \cdots (h_{\epsilon_{1}} -2n+1) \\
& & \cdot e_{\epsilon_1 - \epsilon_2}^{k_{2}}
\cdots e_{\epsilon_1 - \epsilon_{l-1}}^{k_{l-1}}
\cdot (k_{l}! h_{\epsilon_{1}- \epsilon_{l}} \cdots
(h_{\epsilon_{1}- \epsilon_{l}} -k_{l}+1))
f_{\epsilon_1 - \epsilon_{l-1}}^{k_{l-1}}
\cdots f_{\epsilon_1 - \epsilon_2}^{k_{2}}+ U(\frak g){\frak
n}_{+} \\
& \in & (2n)!(-1)^{n-k_{1}}k_{l}!(h_{\epsilon_{1}}-2k_{2}-\ldots
-2k_{l}) \cdots (h_{\epsilon_{1}} -2n+1) \\
& & \cdot (h_{\epsilon_{1}- \epsilon_{l}}-k_{l-1}- \ldots
-k_{2}) \cdots
(h_{\epsilon_{1}- \epsilon_{l}}-k_{l-1}- \ldots
-k_{2} -k_{l}+1) \\
& &  \cdot e_{\epsilon_1 - \epsilon_2}^{k_{2}}
\cdots e_{\epsilon_1 - \epsilon_{l-1}}^{k_{l-1}}
 \cdot  f_{\epsilon_1 - \epsilon_{l-1}}^{k_{l-1}}
\cdots f_{\epsilon_1 - \epsilon_2}^{k_{2}}+ U(\frak g){\frak
n}_{+}  \\
& \ldots \in & (2n)!(-1)^{n-k_{1}}k_{l}! \cdots k_{2}!
(h_{\epsilon_1}-2k_{2}-\ldots -2k_{l})\cdot \ldots \cdot
(h_{\epsilon_1}-2n+1)\cdot \\
& &(h_{\epsilon_1 - \epsilon_l}
-k_{l-1}-\ldots -k_{2}) \cdot \ldots \cdot (h_{\epsilon_1 - \epsilon_l}
-k_{l-1}-\ldots -k_{2}-k_{l}-1) \\
& & \cdot \ldots \cdot
h_{\epsilon_1 - \epsilon_2} \cdot \ldots \cdot (h_{\epsilon_1 -
\epsilon_2}-k_{2}+1) \ + \ U(\frak g){\frak n}_{+}
\end{eqnarray*}
It follows from relation (\ref{2.4.1}) that
\begin{eqnarray*}
& & (e_{\epsilon_1}^{2n})_L
(-\frac{1}{4}f_{\epsilon_1}^2+ f_{\epsilon_1 - \epsilon_2}
f_{\epsilon_1 + \epsilon_2}+ \cdots +f_{\epsilon_1 - \epsilon_l}
f_{\epsilon_1 + \epsilon_l})^{n} \in \\
&  &  (-1)^{n}(2n)!n! {\displaystyle \sum_{{(k_{1},\ldots ,k_{l}) \in {\Z}_{+}^{l} \atop \sum k_{i}=n}}} \frac{1}{k_{1}!4^{k_{1}}} \cdot
(h_{\epsilon_1}-2k_{2}-\ldots -2k_{l})\cdot \ldots \cdot
(h_{\epsilon_1}-2n+1)\cdot \\
& &(h_{\epsilon_1 - \epsilon_l}
-k_{l-1}-\ldots -k_{2}) \cdot \ldots \cdot (h_{\epsilon_1 - \epsilon_l}
-k_{l-1}-\ldots -k_{2}-k_{l}-1) \\
& & \cdot \ldots \cdot
h_{\epsilon_1 - \epsilon_2} \cdot \ldots \cdot (h_{\epsilon_1 -
\epsilon_2}-k_{2}+1) \ + \ U(\frak g){\frak n}_{+}.
\end{eqnarray*}
Finally, relation (\ref{2.4.2}) implies
\begin{eqnarray*}
& &(e_{\epsilon_1}^{2n}f_{\epsilon_1}^{4n})_L
(-\frac{1}{4}e_{\epsilon_1}^2+ e_{\epsilon_1 - \epsilon_2}
e_{\epsilon_1 + \epsilon_2}+ \cdots +e_{\epsilon_1 - \epsilon_l}
e_{\epsilon_1 + \epsilon_l})^n \in \\
&  &  c {\displaystyle \sum_{{(k_{1},\ldots ,k_{l}) \in {\Z}_{+}^{l} \atop \sum k_{i}=n}}} \frac{1}{k_{1}!4^{k_{1}}} \cdot
(h_{\epsilon_1}-2k_{2}-\ldots -2k_{l})\cdot \ldots \cdot
(h_{\epsilon_1}-2n+1)\cdot \\
& &(h_{\epsilon_1 - \epsilon_l}
-k_{l-1}-\ldots -k_{2}) \cdot \ldots \cdot (h_{\epsilon_1 - \epsilon_l}
-k_{l-1}-\ldots -k_{2}-k_{l}-1) \\
& & \cdot \ldots \cdot
h_{\epsilon_1 - \epsilon_2} \cdot \ldots \cdot (h_{\epsilon_1 -
\epsilon_2}-k_{2}+1) \ + \ U(\frak g){\frak n}_{+},
\end{eqnarray*}
for some $c \neq 0$, and the proof is complete.

(2) First notice that
\begin{eqnarray}
&&(f_{\epsilon_1 - \epsilon_2}^{n}f_{\epsilon_1 + \epsilon_2}^{n})_L
(-\frac{1}{4}e_{\epsilon_1}^2+ e_{\epsilon_1 - \epsilon_2}
e_{\epsilon_1 + \epsilon_2}+ \cdots +e_{\epsilon_1 - \epsilon_l}
e_{\epsilon_1 + \epsilon_l})^n \in  R_{0}. \nno
\end{eqnarray}
Lemma \ref{l.2.4.2.1} implies that we can calculate the
corresponding polynomial from
\begin{eqnarray}
&&(-\frac{1}{4}e_{\epsilon_1}^2+ e_{\epsilon_1 - \epsilon_2}
e_{\epsilon_1 + \epsilon_2}+ \cdots +e_{\epsilon_1 - \epsilon_l}
e_{\epsilon_1 + \epsilon_l})^n \, _L
(f_{\epsilon_1 - \epsilon_2}^{n}f_{\epsilon_1 + \epsilon_2}^{n})=
\label{2.4.4} \\
& &= {\displaystyle \sum_{{(k_{1},\ldots ,k_{l}) \in {\Z}_{+}^{l} \atop \sum k_{i}=n}}} {n \choose k_{1},\ldots ,k_{l}}
(-1)^{k_{1}} \frac{1}{4^{k_{1}}} (e_{\epsilon_1}^{2k_{1}}
e_{\epsilon_1 - \epsilon_2}^{k_{2}} e_{\epsilon_1 + \epsilon_2}^{k_{2}}
\cdots e_{\epsilon_1 - \epsilon_l}^{k_{l}} e_{\epsilon_1 +
\epsilon_l}^{k_{l}})_L
(f_{\epsilon_1 - \epsilon_2}^{n}f_{\epsilon_1 + \epsilon_2}^{n}).
\nno
\end{eqnarray}
By using claims (7) and (8) from Lemma \ref{l.2.4.1}, we have
\begin{eqnarray}
&&(e_{\epsilon_1}^{2k_{1}}
e_{\epsilon_1 - \epsilon_2}^{k_{2}} e_{\epsilon_1 + \epsilon_2}^{k_{2}}
\cdots e_{\epsilon_1 - \epsilon_l}^{k_{l}} e_{\epsilon_1 +
\epsilon_l}^{k_{l}})_L (f_{\epsilon_1 - \epsilon_2}^{n}f_{\epsilon_1 +
\epsilon_2}^{n}) = \nno \\
&&=(e_{\epsilon_1 - \epsilon_2}^{k_{2}} \cdots e_{\epsilon_1 - \epsilon_l}^{k_{l}}
e_{\epsilon_1}^{2k_{1}}
e_{\epsilon_1 + \epsilon_2}^{k_{2}} \cdots e_{\epsilon_1 + \epsilon_{l-1}}^{k_{l-1}}
 e_{\epsilon_1 + \epsilon_l}^{k_{l}})_L
(f_{\epsilon_1 + \epsilon_2}^{n}f_{\epsilon_1 - \epsilon_2}^{n}) \in \nno \\
&& \in (e_{\epsilon_1 - \epsilon_2}^{k_{2}} \cdots e_{\epsilon_1 - \epsilon_l}^{k_{l}}
e_{\epsilon_1}^{2k_{1}}
e_{\epsilon_1 + \epsilon_2}^{k_{2}} \cdots  e_{\epsilon_1 + \epsilon_{l-1}}^{k_{l-1}})_L
\Big(n(n-1) \cdots (n-k_{l}+1) \cdot \nno \\
&& \qquad  \cdot f_{\epsilon_1 + \epsilon_2}^{n-k_{l}}f_{\epsilon_2 - \epsilon_{l}}^{k_{l}}
f_{\epsilon_1 - \epsilon_2}^{n} \Big)  +  U(\frak g){\frak n}_{+} = \nno \\
&& = (e_{\epsilon_1 - \epsilon_2}^{k_{2}} \cdots e_{\epsilon_1 - \epsilon_l}^{k_{l}}
e_{\epsilon_1}^{2k_{1}}
e_{\epsilon_1 + \epsilon_2}^{k_{2}} \cdots  e_{\epsilon_1 + \epsilon_{l-2}}^{k_{l-2}})_L
\Big(n(n-1) \cdots (n-k_{l}+1) \cdot \nno \\
&& \qquad  \cdot (n-k_{l}) \cdots (n-k_{l}-k_{l-1}+1) f_{\epsilon_1 + \epsilon_2}^{n-k_{l}-k_{l-1}}
f_{\epsilon_2 - \epsilon_{l-1}}^{k_{l-1}}
f_{\epsilon_2 - \epsilon_{l}}^{k_{l}}
f_{\epsilon_1 - \epsilon_2}^{n}\Big)  +  U(\frak g){\frak n}_{+} = \nno \\
&& \ldots = (e_{\epsilon_1 - \epsilon_2}^{k_{2}} \cdots e_{\epsilon_1 - \epsilon_l}^{k_{l}}
e_{\epsilon_1}^{2k_{1}} e_{\epsilon_1 + \epsilon_2}^{k_{2}})_L
\Big(n(n-1) \cdots (n-k_{l}- \ldots -k_{3}+1) \cdot \nno \\
&& \qquad  \cdot  f_{\epsilon_1 + \epsilon_2}^{n-k_{l}- \ldots -k_{3}}
f_{\epsilon_2 - \epsilon_{3}}^{k_{3}} \cdots
f_{\epsilon_2 - \epsilon_{l-1}}^{k_{l-1}}f_{\epsilon_2 - \epsilon_{l}}^{k_{l}}
f_{\epsilon_1 - \epsilon_2}^{n}\Big)  +  U(\frak g){\frak n}_{+} = \nno \\
&& = (e_{\epsilon_1 - \epsilon_2}^{k_{2}} \cdots e_{\epsilon_1 - \epsilon_l}^{k_{l}}
e_{\epsilon_1}^{2k_{1}} e_{\epsilon_1 + \epsilon_2}^{k_{2}})_L
\Big(n(n-1) \cdots (k_{1}+k_{2}+1) \cdot \nno \\
&& \qquad  \cdot  f_{\epsilon_1 + \epsilon_2}^{k_{1}+k_{2}}
f_{\epsilon_2 - \epsilon_{3}}^{k_{3}} \cdots
f_{\epsilon_2 - \epsilon_{l-1}}^{k_{l-1}}f_{\epsilon_2 - \epsilon_{l}}^{k_{l}}
f_{\epsilon_1 - \epsilon_2}^{n}\Big)  +  U(\frak g){\frak n}_{+}. \nno
\end{eqnarray}
Claims (9), (8), (10) and (6) from Lemma \ref{l.2.4.1} imply
\begin{eqnarray}
&&(e_{\epsilon_1 - \epsilon_2}^{k_{2}} \cdots e_{\epsilon_1 - \epsilon_l}^{k_{l}}
e_{\epsilon_1}^{2k_{1}}
e_{\epsilon_1 + \epsilon_2}^{k_{2}} \cdots e_{\epsilon_1 + \epsilon_{l-1}}^{k_{l-1}}
 e_{\epsilon_1 + \epsilon_l}^{k_{l}})_L
(f_{\epsilon_1 + \epsilon_2}^{n}f_{\epsilon_1 - \epsilon_2}^{n}) \in \nno \\
&& \in (e_{\epsilon_1 - \epsilon_2}^{k_{2}} \cdots e_{\epsilon_1 - \epsilon_l}^{k_{l}}
e_{\epsilon_1}^{2k_{1}})_L
\Big(n(n-1) \cdots (k_{1}+ k_{2}+1) \cdot \nno \\
&& \qquad  \cdot (e_{\epsilon_1 + \epsilon_2}^{k_{2}})_L( f_{\epsilon_1 + \epsilon_2}^{k_{1}+k_{2}})
\cdot f_{\epsilon_2 - \epsilon_{3}}^{k_{3}} \cdots
f_{\epsilon_2 - \epsilon_{l-1}}^{k_{l-1}}f_{\epsilon_2 - \epsilon_{l}}^{k_{l}}
f_{\epsilon_1 - \epsilon_2}^{n}\Big)  +  U(\frak g){\frak n}_{+} = \nno \\
&& = (e_{\epsilon_1 - \epsilon_2}^{k_{2}} \cdots e_{\epsilon_1 - \epsilon_l}^{k_{l}}
e_{\epsilon_1}^{2k_{1}})_L
\Big(n \cdots (k_{1}+ k_{2}+1) \cdot (k_{1}+ k_{2}) \cdots (k_{1}+1)
f_{\epsilon_1 + \epsilon_2}^{k_{1}} \cdot  \nno \\
&& \quad  \cdot (h_{\epsilon_1 +
\epsilon_2}-k_{1}) \cdots (h_{\epsilon_1 + \epsilon_2}-k_{1}-k_{2}+1)
f_{\epsilon_2 - \epsilon_{3}}^{k_{3}} \cdots
f_{\epsilon_2 - \epsilon_{l-1}}^{k_{l-1}}f_{\epsilon_2 - \epsilon_{l}}^{k_{l}}
f_{\epsilon_1 - \epsilon_2}^{n}\Big)  +  U(\frak g){\frak n}_{+}  \nno \\
&& = (e_{\epsilon_1 - \epsilon_2}^{k_{2}} \cdots e_{\epsilon_1 - \epsilon_l}^{k_{l}}
e_{\epsilon_1}^{2k_{1}})_L
\Big(n(n-1) \cdots (k_{1}+1) f_{\epsilon_1 + \epsilon_2}^{k_{1}}
f_{\epsilon_2 - \epsilon_{3}}^{k_{3}} \cdots
f_{\epsilon_2 - \epsilon_{l-1}}^{k_{l-1}}f_{\epsilon_2 - \epsilon_{l}}^{k_{l}}
f_{\epsilon_1 - \epsilon_2}^{n} \cdot  \nno \\
&& \quad  \cdot (h_{\epsilon_1 +
\epsilon_2}-k_{1}-k_{3}- \ldots - k_{l}) \cdots (h_{\epsilon_1 +
\epsilon_2}-k_{1}-k_{2}-k_{3}- \ldots - k_{l}+1) \Big)  +  U(\frak g){\frak n}_{+}   \nno \\
&& = (e_{\epsilon_1 - \epsilon_2}^{k_{2}} \cdots e_{\epsilon_1 - \epsilon_l}^{k_{l}}
e_{\epsilon_1}^{2k_{1}})_L
\Big(n(n-1) \cdots (k_{1}+1) f_{\epsilon_1 + \epsilon_2}^{k_{1}}
f_{\epsilon_2 - \epsilon_{3}}^{k_{3}} \cdots
f_{\epsilon_2 - \epsilon_{l-1}}^{k_{l-1}}f_{\epsilon_2 - \epsilon_{l}}^{k_{l}}
f_{\epsilon_1 - \epsilon_2}^{n} \cdot \nno \\
&& \quad  \cdot (h_{\epsilon_1 +
\epsilon_2}-n+ k_{2}) \cdots (h_{\epsilon_1 +
\epsilon_2}-n+1) \Big)  +  U(\frak g){\frak
n}_{+}. \nno
\end{eqnarray}
It follows from claims (3) and (5) from Lemma \ref{l.2.4.1} that
\begin{eqnarray}
&&(e_{\epsilon_1 - \epsilon_2}^{k_{2}} \cdots e_{\epsilon_1 - \epsilon_l}^{k_{l}}
e_{\epsilon_1}^{2k_{1}}
e_{\epsilon_1 + \epsilon_2}^{k_{2}} \cdots e_{\epsilon_1 + \epsilon_{l-1}}^{k_{l-1}}
 e_{\epsilon_1 + \epsilon_l}^{k_{l}})_L
(f_{\epsilon_1 + \epsilon_2}^{n}f_{\epsilon_1 - \epsilon_2}^{n}) \in \nno \\
&& \in (e_{\epsilon_1 - \epsilon_2}^{k_{2}} \cdots e_{\epsilon_1 - \epsilon_l}^{k_{l}})_L
\Big(n(n-1) \cdots (k_{1}+1)(-1)^{k_{1}}(2k_{1})! \, e_{\epsilon_1 - \epsilon_2}^{k_{1}}
f_{\epsilon_2 - \epsilon_{3}}^{k_{3}} \cdots f_{\epsilon_2 - \epsilon_{l}}^{k_{l}}
f_{\epsilon_1 - \epsilon_2}^{n} \cdot  \nno \\
&& \quad  \cdot (h_{\epsilon_1 + \epsilon_2}-n+ k_{2}) \cdots (h_{\epsilon_1 +
\epsilon_2}-n+1) \Big)  +  U(\frak g){\frak n}_{+}. \nno
\end{eqnarray}
By using claims (11), (12), (10) and (1) from Lemma \ref{l.2.4.1}, we obtain
\begin{eqnarray}
&&(e_{\epsilon_1 - \epsilon_2}^{k_{2}} \cdots e_{\epsilon_1 - \epsilon_l}^{k_{l}}
e_{\epsilon_1}^{2k_{1}}
e_{\epsilon_1 + \epsilon_2}^{k_{2}} \cdots e_{\epsilon_1 + \epsilon_{l-1}}^{k_{l-1}}
 e_{\epsilon_1 + \epsilon_l}^{k_{l}})_L
(f_{\epsilon_1 + \epsilon_2}^{n}f_{\epsilon_1 - \epsilon_2}^{n}) \in \nno \\
&& \in (e_{\epsilon_1 - \epsilon_2}^{k_{2}} \cdots e_{\epsilon_1 - \epsilon_{l-1}}^{k_{l-1}})_L
\Big(n(n-1) \cdots (k_{1}+1)(-1)^{k_{1}}(2k_{1})! \, e_{\epsilon_1 - \epsilon_2}^{k_{1}}
f_{\epsilon_2 - \epsilon_{3}}^{k_{3}} \cdots f_{\epsilon_2 - \epsilon_{l-1}}^{k_{l-1}}
\cdot  \nno \\
&& \quad  \cdot k_{l}! \, e_{\epsilon_1 - \epsilon_2}^{k_{l}}
f_{\epsilon_1 - \epsilon_2}^{n} (h_{\epsilon_1 + \epsilon_2}-n+ k_{2}) \cdots (h_{\epsilon_1 +
\epsilon_2}-n+1) \Big)  +  U(\frak g){\frak n}_{+} = \nno \\
&& = (e_{\epsilon_1 - \epsilon_2}^{k_{2}} \cdots e_{\epsilon_1 - \epsilon_{l-2}}^{k_{l-2}})_L
\Big(n(n-1) \cdots (k_{1}+1)(-1)^{k_{1}}(2k_{1})! \, e_{\epsilon_1 - \epsilon_2}^{k_{1}}
f_{\epsilon_2 - \epsilon_{3}}^{k_{3}} \cdots f_{\epsilon_2 - \epsilon_{l-2}}^{k_{l-2}}
\cdot  \nno \\
&& \quad  \cdot k_{l-1}! \, e_{\epsilon_1 - \epsilon_2}^{k_{l-1}}
k_{l}! \, e_{\epsilon_1 - \epsilon_2}^{k_{l}}
f_{\epsilon_1 - \epsilon_2}^{n} (h_{\epsilon_1 + \epsilon_2}-n+ k_{2}) \cdots (h_{\epsilon_1 +
\epsilon_2}-n+1) \Big)  +  U(\frak g){\frak n}_{+} = \nno \\
&& \ldots = (e_{\epsilon_1 - \epsilon_2}^{k_{2}})_L
\Big(n(n-1) \cdots (k_{1}+1)(-1)^{k_{1}}(2k_{1})! \,k_{3}! \cdots k_{l}! \,
e_{\epsilon_1 - \epsilon_2}^{k_{1}}e_{\epsilon_1 -
\epsilon_2}^{k_{3}} \ldots
e_{\epsilon_1 - \epsilon_2}^{k_{l}} \cdot  \nno \\
&& \quad  \cdot f_{\epsilon_1 - \epsilon_2}^{n}
(h_{\epsilon_1 + \epsilon_2}-n+ k_{2}) \cdots (h_{\epsilon_1 +
\epsilon_2}-n+1) \Big)  +  U(\frak g){\frak n}_{+} = \nno \\
&&  = n(n-1) \cdots (k_{1}+1)(-1)^{k_{1}}(2k_{1})! \, k_{3}! \cdots k_{l}! \,
e_{\epsilon_1 - \epsilon_2}^{n-k_{2}} (e_{\epsilon_1 - \epsilon_2}^{k_{2}})_L (f_{\epsilon_1 - \epsilon_2}^{n})
 \cdot  \nno \\
&& \quad   \cdot
(h_{\epsilon_1 + \epsilon_2}-n+ k_{2}) \cdots (h_{\epsilon_1 +
\epsilon_2}-n+1)   +  U(\frak g){\frak n}_{+} = \nno \\
&&  = n(n-1) \cdots (k_{1}+1)(-1)^{k_{1}}(2k_{1})! \, k_{3}! \cdots k_{l}! \,
e_{\epsilon_1 - \epsilon_2}^{n-k_{2}} n(n-1) \cdots (n-k_{2}+1)
f_{\epsilon_1 - \epsilon_2}^{n-k_{2}} \cdot   \nno \\
&& \quad   \cdot (h_{\epsilon_1 - \epsilon_2}-n+ k_{2})\cdots (h_{\epsilon_1
- \epsilon_2}-n+1)
(h_{\epsilon_1 + \epsilon_2}-n+ k_{2}) \cdots (h_{\epsilon_1 +
\epsilon_2}-n+1)   +  U(\frak g){\frak n}_{+}  \nno \\
&& = \frac{n!}{k_{1}!}(-1)^{k_{1}}(2k_{1})! \, k_{3}! \cdots k_{l}! \,
n! \, h_{\epsilon_1 - \epsilon_2} \cdots  (h_{\epsilon_1 - \epsilon_2}-n+1)
\cdot \nno \\
&& \qquad \cdot (h_{\epsilon_1 + \epsilon_2}-n+ k_{2}) \cdots (h_{\epsilon_1 +
\epsilon_2}-n+1)   +  U(\frak g){\frak n}_{+}.  \nno
\end{eqnarray}
It follows from relation (\ref{2.4.4}) that
\begin{eqnarray}
&&(-\frac{1}{4}e_{\epsilon_1}^2+ e_{\epsilon_1 - \epsilon_2}
e_{\epsilon_1 + \epsilon_2}+ \cdots +e_{\epsilon_1 - \epsilon_l}
e_{\epsilon_1 + \epsilon_l})^n \, _L
(f_{\epsilon_1 - \epsilon_2}^{n}f_{\epsilon_1 + \epsilon_2}^{n}) \in  \nno \\
& &\in (n!)^{3} h_{\epsilon_1 - \epsilon_2} \cdots  (h_{\epsilon_1 - \epsilon_2}-n+1)
\cdot  \nno \\
&& \quad \cdot {\displaystyle \sum_{{(k_{1},\ldots ,k_{l}) \in {\Z}_{+}^{l} \atop \sum k_{i}=n}}}
{2k_{1} \choose k_{1}} \frac{1}{4^{k_{1}}}\frac{1}{k_{2}!}
(h_{\epsilon_1 + \epsilon_2}-n+ k_{2}) \cdots (h_{\epsilon_1 + \epsilon_2}-n+1)
 +   U(\frak g){\frak n}_{+}. \nno
\end{eqnarray}
Since
\begin{eqnarray*}
&& {\displaystyle \sum_{{(k_{1},\ldots ,k_{l}) \in {\Z}_{+}^{l} \atop \sum k_{i}=n}}}
{2k_{1} \choose k_{1}} \frac{1}{4^{k_{1}}}\frac{1}{k_{2}!}
(h_{\epsilon_1 + \epsilon_2}-n+ k_{2}) \cdots (h_{\epsilon_1 +
\epsilon_2}-n+1)= \\
&& ={\displaystyle \sum_{{(k_{1},k_{2}) \in {\Z}_{+}^{2} \atop  k_{1}+k_{2} \leq n}}}
{n-k_{1}-k_{2}+l-3 \choose n-k_{1}-k_{2}}{2k_{1} \choose k_{1}} \frac{1}{4^{k_{1}}}
{h_{\epsilon_1 + \epsilon_2}-n+ k_{2} \choose k_{2}}= \\
&& ={\displaystyle \sum_{k_{2}=0}^{n}} {h_{\epsilon_1 + \epsilon_2}-n+ k_{2} \choose k_{2}}
{\displaystyle \sum_{k_{1}=0}^{n-k_{2}}}
{n-k_{1}-k_{2}+l-3 \choose n-k_{1}-k_{2}}{2k_{1} \choose k_{1}} \frac{1}{4^{k_{1}}}= \\
&& ={\displaystyle \sum_{k_{2}=0}^{n}}  {h_{\epsilon_1 + \epsilon_2}-n+ k_{2} \choose k_{2}}
{n-k_{2}+l-\frac{5}{2}  \choose n-k_{2}}= {h_{\epsilon_1 + \epsilon_2}+l- \frac{3}{2}
\choose n},
\end{eqnarray*}
we get
\begin{eqnarray}
&&(f_{\epsilon_1 - \epsilon_2}^{n}f_{\epsilon_1 + \epsilon_2}^{n})_L
(-\frac{1}{4}e_{\epsilon_1}^2+ e_{\epsilon_1 - \epsilon_2}
e_{\epsilon_1 + \epsilon_2}+ \cdots +e_{\epsilon_1 - \epsilon_l}
e_{\epsilon_1 + \epsilon_l})^n \in \nno \\
&& \in c_{1} h_{\epsilon_1 - \epsilon_2}
\cdot \ldots \cdot (h_{\epsilon_1 - \epsilon_2}-n+1)
(h_{\epsilon_1 + \epsilon_{2}}+l-\frac{3}{2})
\cdot \ldots
\cdot (h_{\epsilon_1 + \epsilon_{2}}+l-n- \frac{1}{2})
+ U(\frak g){\frak n}_{+}  \nno \\
&& =  c_{1} \, p_{1}(h)  + U(\frak g){\frak n}_{+},
\label{2.4.5}
\end{eqnarray}
for some $c_{1} \neq 0$.

Thus, $p_{1} \in {\mathcal P}_{0}$. Let's prove $p_{i} \in {\mathcal
P}_{0}$, for $1=2, \ldots ,l-1$.  Using notation
$$ \bar{u}= -\frac{1}{4}e_{\epsilon_1}^2+ e_{\epsilon_1 - \epsilon_2}
e_{\epsilon_1 + \epsilon_2}+ \cdots +e_{\epsilon_1 - \epsilon_l}
e_{\epsilon_1 + \epsilon_l}, $$
as in (1), we have
$$(f_{\epsilon_1 +\epsilon_{i}}^{2})_{L} \bar{u}=
2(-\frac{1}{4}f_{\epsilon_{i}}^2- e_{\epsilon_1 - \epsilon_{i}}
f_{\epsilon_1 + \epsilon_{i}}- \cdots -e_{\epsilon_{i-1} - \epsilon_{i}}
f_{\epsilon_{i-1} + \epsilon_{i}}+
f_{\epsilon_{i} - \epsilon_{i+1}}
f_{\epsilon_{i} + \epsilon_{i+1}}+ \cdots + f_{\epsilon_{i} - \epsilon_{l}}
f_{\epsilon_{i} + \epsilon_{l}} )$$
and $(f_{\epsilon_1 +\epsilon_{i}}^{3})_{L} \bar{u}=0,$
for $1=2, \ldots ,l-1$, which implies
\begin{eqnarray*}
& &(f_{\epsilon_1 +\epsilon_{i}}^{2n})_L
(-\frac{1}{4}e_{\epsilon_1}^2+ e_{\epsilon_1 - \epsilon_2}
e_{\epsilon_1 + \epsilon_2}+ \cdots +e_{\epsilon_1 - \epsilon_l}
e_{\epsilon_1 + \epsilon_l})^n=
(f_{\epsilon_1 +\epsilon_{i}}^{2n})_L(\bar{u}^{n}) =  \\
& &=\frac{(2n)!}{(2!)^{n}}
((f_{\epsilon_1 +\epsilon_{i}}^{2})_L\bar{u})^{n} =
(2n)! (-\frac{1}{4}f_{\epsilon_{i}}^2- e_{\epsilon_1 - \epsilon_{i}}
f_{\epsilon_1 + \epsilon_{i}}- \cdots -e_{\epsilon_{i-1} - \epsilon_{i}}
f_{\epsilon_{i-1} + \epsilon_{i}}+ \\
&& \quad + f_{\epsilon_{i} - \epsilon_{i+1}}
f_{\epsilon_{i} + \epsilon_{i+1}}+ \cdots + f_{\epsilon_{i} - \epsilon_{l}}
f_{\epsilon_{i} + \epsilon_{l}} )^{n}.
\end{eqnarray*}
Since all root vectors $f_{\epsilon_i},
e_{\epsilon_1 - \epsilon_{i}}, f_{\epsilon_1 + \epsilon_{i}}, \ldots ,
e_{\epsilon_{i-1} - \epsilon_{i}}$, $f_{\epsilon_{i-1} + \epsilon_{i}},
f_{\epsilon_{i} - \epsilon_{i+1}}, f_{\epsilon_{i} +
\epsilon_{i+1}},$ $\ldots , f_{\epsilon_{i} - \epsilon_{l}},
f_{\epsilon_{i} + \epsilon_{l}}  \in \frak g$
commute, we get
\begin{eqnarray}
& &(-\frac{1}{4}f_{\epsilon_{i}}^2- e_{\epsilon_1 - \epsilon_{i}}
f_{\epsilon_1 + \epsilon_{i}}- \cdots -e_{\epsilon_{i-1} - \epsilon_{i}}
f_{\epsilon_{i-1} + \epsilon_{i}}+
f_{\epsilon_{i} - \epsilon_{i+1}}
f_{\epsilon_{i} + \epsilon_{i+1}}+ \cdots + f_{\epsilon_{i} - \epsilon_{l}}
f_{\epsilon_{i} + \epsilon_{l}} )^{n}= \nonumber \\
& &= {\displaystyle \sum_{{(k_{1},\ldots ,k_{l}) \in {\Z}_{+}^{l} \atop \sum k_{j}=n}}} {n \choose k_{1},\ldots ,k_{l}}
(-1)^{k_{1}+ \ldots +k_{i}} \frac{1}{4^{k_{i}}}
f_{\epsilon_{i}}^{2k_{i}}
e_{\epsilon_1 - \epsilon_{i}}^{k_{1}} f_{\epsilon_1 + \epsilon_{i}}^{k_{1}}
\cdots e_{\epsilon_{i-1} - \epsilon_{i}}^{k_{i-1}} f_{\epsilon_{i-1}
+ \epsilon_{i}}^{k_{i-1}} \cdot \nno \\
&& \qquad \qquad \qquad \cdot f_{\epsilon_{i} - \epsilon_{i+1}}^{k_{i+1}} f_{\epsilon_{i} + \epsilon_{i+1}}^{k_{i+1}}
\cdots f_{\epsilon_{i} - \epsilon_{l}}^{k_{l}}f_{\epsilon_{i} +
\epsilon_{l}}^{k_{l}}. \label{2.4.6}
\end{eqnarray}
If any of indices $k_{1}, \ldots , k_{i-1}$  is nonzero,
then the corresponding summand in (\ref{2.4.6})
is an element of $U(\frak g){\frak n}_{+}$. Thus we obtain
\begin{eqnarray}
& &(f_{\epsilon_1 +\epsilon_{i}}^{2n})_L
(-\frac{1}{4}e_{\epsilon_1}^2+ e_{\epsilon_1 - \epsilon_2}
e_{\epsilon_1 + \epsilon_2}+ \cdots +e_{\epsilon_1 - \epsilon_l}
e_{\epsilon_1 + \epsilon_l})^n \in \nno \\
& & \in (2n)! \, {\displaystyle \sum_{{(k_{i},\ldots ,k_{l}) \in {\Z}_{+}^{l} \atop \sum k_{j}=n}}} {n \choose k_{i},\ldots ,k_{l}}
(-1)^{k_{i}} \frac{1}{4^{k_{i}}}
f_{\epsilon_{i}}^{2k_{i}}
f_{\epsilon_{i} - \epsilon_{i+1}}^{k_{i+1}} f_{\epsilon_{i} + \epsilon_{i+1}}^{k_{i+1}}
\cdots f_{\epsilon_{i} - \epsilon_{l}}^{k_{l}}f_{\epsilon_{i} +
\epsilon_{l}}^{k_{l}} + U(\frak g){\frak n}_{+} \nno \\
& & = (2n)! \,(-\frac{1}{4}f_{\epsilon_{i}}^2+
f_{\epsilon_{i} - \epsilon_{i+1}}
f_{\epsilon_{i} + \epsilon_{i+1}}+ \cdots + f_{\epsilon_{i} - \epsilon_{l}}
f_{\epsilon_{i} + \epsilon_{l}} )^{n} + U(\frak g){\frak n}_{+}. \nonumber
\end{eqnarray}

Let $\frak g'$ be the subalgebra of $\frak g$ associated to roots
$\alpha_{i}, \ldots ,\alpha_{l}$. Then $\frak g'$ is a simple Lie algebra
of type $B_{l-i+1}$. Let ${\frak g}'={\frak n}_{-}' \oplus {\frak h}'
\oplus {\frak n}_{+}'$ be the corresponding triangular decomposition
of ${\frak g}'$. Universal enveloping algebra $U(\frak g')$
is embedded in $U(\frak g)$ in the natural way.
Vector
\begin{eqnarray}
(-\frac{1}{4}f_{\epsilon_{i}}^2+
f_{\epsilon_{i} - \epsilon_{i+1}}
f_{\epsilon_{i} + \epsilon_{i+1}}+ \cdots + f_{\epsilon_{i} - \epsilon_{l}}
f_{\epsilon_{i} + \epsilon_{l}} )^{n} \in U(\frak g') \nonumber
\end{eqnarray}
is the lowest weight vector for $\frak g'$. Let $R'$ be a $\frak
g'$-module generated with this vector, and $R_{0}'$ zero-weight
subspace of $R'$. $R'$ is then a $\frak g'$-module with highest weight
$2 \epsilon_{i}$ and highest weight vector
\begin{eqnarray}
(-\frac{1}{4}e_{\epsilon_{i}}^2+
e_{\epsilon_{i} - \epsilon_{i+1}}
e_{\epsilon_{i} + \epsilon_{i+1}}+ \cdots + e_{\epsilon_{i} - \epsilon_{l}}
e_{\epsilon_{i} + \epsilon_{l}} )^{n}. \nonumber
\end{eqnarray}
Clearly,
\begin{eqnarray}
&&(f_{\epsilon_{i} - \epsilon_{i+1}}^{n}f_{\epsilon_{i} + \epsilon_{i+1}}^{n})_L
(-\frac{1}{4}e_{\epsilon_{i}}^2+
e_{\epsilon_{i} - \epsilon_{i+1}}
e_{\epsilon_{i} + \epsilon_{i+1}}+ \cdots + e_{\epsilon_{i} - \epsilon_{l}}
e_{\epsilon_{i} + \epsilon_{l}} )^{n} \in  R_{0}'. \nno
\end{eqnarray}
If we apply relation (\ref{2.4.5}) to subalgebra $\frak g'$ of
type $B_{l-i+1}$, we get
\begin{eqnarray}
&&(f_{\epsilon_{i} - \epsilon_{i+1}}^{n}f_{\epsilon_{i} + \epsilon_{i+1}}^{n})_L
(-\frac{1}{4}e_{\epsilon_{i}}^2+
e_{\epsilon_{i} - \epsilon_{i+1}}
e_{\epsilon_{i} + \epsilon_{i+1}}+ \cdots + e_{\epsilon_{i} - \epsilon_{l}}
e_{\epsilon_{i} + \epsilon_{l}} )^{n} \in \nno \\
&&  \in c_{1} h_{\epsilon_{i} - \epsilon_{i+1}}
\cdot \ldots \cdot (h_{\epsilon_{i} - \epsilon_{i+1}}-n+1)
(h_{\epsilon_{i} + \epsilon_{i+1}}+l-i-\frac{1}{2})
\cdot \ldots \cdot \nno \\
&& \qquad \cdot (h_{\epsilon_{i} + \epsilon_{i+1}}+l-n-i+ \frac{1}{2}) + U(\frak g'){\frak n'}_{+}
=  c_{1} \, p_{i}(h)  + U(\frak g'){\frak n'}_{+}. \nno
\end{eqnarray}
Clearly, there exists $Y \in U({\frak n'}_{+})$ such that
\begin{eqnarray}
&&Y_L
(-\frac{1}{4}f_{\epsilon_{i}}^2+
f_{\epsilon_{i} - \epsilon_{i+1}}
f_{\epsilon_{i} + \epsilon_{i+1}}+ \cdots + f_{\epsilon_{i} - \epsilon_{l}}
f_{\epsilon_{i} + \epsilon_{l}} )^{n}= \nno \\
&&=(f_{\epsilon_{i} - \epsilon_{i+1}}^{n}f_{\epsilon_{i} + \epsilon_{i+1}}^{n})_L
(-\frac{1}{4}e_{\epsilon_{i}}^2+
e_{\epsilon_{i} - \epsilon_{i+1}}
e_{\epsilon_{i} + \epsilon_{i+1}}+ \cdots + e_{\epsilon_{i} - \epsilon_{l}}
e_{\epsilon_{i} + \epsilon_{l}} )^{n},\nno
\end{eqnarray}
which implies
\begin{eqnarray}
& &(Y f_{\epsilon_1 +\epsilon_{i}}^{2n})_L
(-\frac{1}{4}e_{\epsilon_1}^2+ e_{\epsilon_1 - \epsilon_2}
e_{\epsilon_1 + \epsilon_2}+ \cdots +e_{\epsilon_1 - \epsilon_l}
e_{\epsilon_1 + \epsilon_l})^n \in \nno \\
& & \in (2n)! \,Y_L (-\frac{1}{4}f_{\epsilon_{i}}^2+
f_{\epsilon_{i} - \epsilon_{i+1}}
f_{\epsilon_{i} + \epsilon_{i+1}}+ \cdots + f_{\epsilon_{i} - \epsilon_{l}}
f_{\epsilon_{i} + \epsilon_{l}} )^{n} + U(\frak g){\frak n}_{+}
\nonumber \\
&& = c_{i} \, p_{i}(h)  + U(\frak g'){\frak n'}_{+} + U(\frak g){\frak n}_{+}
=  c_{i} \, p_{i}(h)  + U(\frak g){\frak n}_{+},
\end{eqnarray}
for some $c_{i} \neq 0$. Thus $p_{i} \in {\mathcal P}_{0}$,
for $i=2, \ldots ,l-1$.

(3) We claim that
\begin{eqnarray*}
&&(e_{\epsilon_l}^{2n}f_{\epsilon_1 + \epsilon_l}^{2n})_L
(-\frac{1}{4}e_{\epsilon_1}^2+ e_{\epsilon_1 - \epsilon_2}
e_{\epsilon_1 + \epsilon_2}+ \cdots +e_{\epsilon_1 - \epsilon_l}
e_{\epsilon_1 + \epsilon_l})^n \in  c_{l} \, p_{l}(h)  +  U(\frak g){\frak
n}_{+},
\end{eqnarray*}
for some $c_{l} \neq 0$.

Using notation
$$ \bar{u}= -\frac{1}{4}e_{\epsilon_1}^2+ e_{\epsilon_1 - \epsilon_2}
e_{\epsilon_1 + \epsilon_2}+ \cdots +e_{\epsilon_1 - \epsilon_l}
e_{\epsilon_1 + \epsilon_l} \in U(\frak g), $$
as in $(1)$, we have
$$(f_{\epsilon_1 + \epsilon_l}^{2})_{L} \bar{u}=
-2(\frac{1}{4}f_{\epsilon_l}^2+ f_{\epsilon_1 + \epsilon_l}
e_{\epsilon_1 - \epsilon_l}+ f_{\epsilon_2 + \epsilon_l}
e_{\epsilon_2 - \epsilon_l}+ \cdots +f_{\epsilon_{l-1} + \epsilon_l}
e_{\epsilon_{l-1} - \epsilon_l})$$
and $(f_{\epsilon_1 + \epsilon_l}^{3})_{L} \bar{u}=0$,
which implies
\begin{eqnarray*}
& &(f_{\epsilon_1 + \epsilon_l}^{2n})_L
(-\frac{1}{4}e_{\epsilon_1}^2+ e_{\epsilon_1 - \epsilon_2}
e_{\epsilon_1 + \epsilon_2}+ \cdots +e_{\epsilon_1 - \epsilon_l}
e_{\epsilon_1 + \epsilon_l})^n=
(f_{\epsilon_1 + \epsilon_l}^{2n})_L(\bar{u}^{n})= \\
& &= \frac{(2n)!}{(2!)^{n}}
((f_{\epsilon_1 + \epsilon_l}^{2})_L\bar{u})^{n}
 =(-1)^{n}(2n)! (\frac{1}{4}f_{\epsilon_l}^2+ f_{\epsilon_1 + \epsilon_l}
e_{\epsilon_1 - \epsilon_l}+ \cdots +f_{\epsilon_{l-1} + \epsilon_l}
e_{\epsilon_{l-1} - \epsilon_l})^{n}.
\end{eqnarray*}
We obtain
\begin{eqnarray*}
& &(e_{\epsilon_l}^{2n}f_{\epsilon_1 + \epsilon_l}^{2n})_L
(-\frac{1}{4}e_{\epsilon_1}^2+ e_{\epsilon_1 - \epsilon_2}
e_{\epsilon_1 + \epsilon_2}+ \cdots +e_{\epsilon_1 - \epsilon_l}
e_{\epsilon_1 + \epsilon_l})^n = \nonumber \\
& &\quad  = (-1)^{n}(2n)!(e_{\epsilon_l}^{2n})_L
(\frac{1}{4}f_{\epsilon_l}^2+ f_{\epsilon_1 + \epsilon_l}
e_{\epsilon_1 - \epsilon_l}+ f_{\epsilon_2 + \epsilon_l}
e_{\epsilon_2 - \epsilon_l}+ \cdots +f_{\epsilon_{l-1} + \epsilon_l}
e_{\epsilon_{l-1} - \epsilon_l})^{n}.
\end{eqnarray*}
Clearly,
\begin{eqnarray}
& &(\frac{1}{4}f_{\epsilon_l}^2+ f_{\epsilon_1 + \epsilon_l}
e_{\epsilon_1 - \epsilon_l}+ f_{\epsilon_2 + \epsilon_l}
e_{\epsilon_2 - \epsilon_l}+ \cdots +f_{\epsilon_{l-1} + \epsilon_l}
e_{\epsilon_{l-1} - \epsilon_l})^{n}=  \\
& & \ = {\displaystyle \sum_{{(k_{1},\ldots ,k_{l}) \in {\Z}_{+}^{l} \atop \sum k_{i}=n}}} {n \choose k_{1},\ldots ,k_{l}}
\frac{1}{4^{k_{l}}} f_{\epsilon_l}^{2k_{l}}
f_{\epsilon_1 + \epsilon_l}^{k_{1}} e_{\epsilon_1 - \epsilon_l}^{k_{1}}
\cdots f_{\epsilon_{l-1} + \epsilon_l}^{k_{l-1}} e_{\epsilon_{l-1} -
\epsilon_l}^{k_{l-1}}= \nonumber \\
& & \ = {\displaystyle \sum_{{(k_{1},\ldots ,k_{l}) \in {\Z}_{+}^{l} \atop \sum k_{i}=n}}} {n \choose k_{1},\ldots ,k_{l}}
\frac{1}{4^{k_{l}}} f_{\epsilon_l}^{2k_{l}}
f_{\epsilon_1 + \epsilon_l}^{k_{1}} \cdots f_{\epsilon_{l-1} + \epsilon_l}^{k_{l-1}}
e_{\epsilon_1 - \epsilon_l}^{k_{1}} \cdots e_{\epsilon_{l-1} -
\epsilon_l}^{k_{l-1}} \nonumber
\end{eqnarray}
If any of indices $k_{1} \ldots k_{l-1}$ is nonzero,
then the corresponding summand above is an element of
$U(\frak g){\frak n}_{+}$. It follows that
\begin{eqnarray*}
& & (e_{\epsilon_l}^{2n}f_{\epsilon_1 + \epsilon_l}^{2n})_L
(-\frac{1}{4}e_{\epsilon_1}^2+ e_{\epsilon_1 - \epsilon_2}
e_{\epsilon_1 + \epsilon_2}+ \cdots +e_{\epsilon_1 - \epsilon_l}
e_{\epsilon_1 + \epsilon_l})^n \in \\
& & \qquad  \in (-1)^{n}(2n)! \frac{1}{4^{n}} (e_{\epsilon_l}^{2n})_L
(f_{\epsilon_l}^{2n}) + U(\frak g){\frak n}_{+}.
\end{eqnarray*}
By using claim (1) from Lemma \ref{l.2.4.1}, we obtain
\begin{eqnarray*}
& & (e_{\epsilon_l}^{2n}f_{\epsilon_1 + \epsilon_l}^{2n})_L
(-\frac{1}{4}e_{\epsilon_1}^2+ e_{\epsilon_1 - \epsilon_2}
e_{\epsilon_1 + \epsilon_2}+ \cdots +e_{\epsilon_1 - \epsilon_l}
e_{\epsilon_1 + \epsilon_l})^n \in \\
& & \qquad  \in  c_{l} \, h_{l}(h_{l}-1)\cdot \ldots \cdot (h_{l}-2n+1)
+  U(\frak g){\frak n}_{+},
\end{eqnarray*}
for some $c_{l} \neq 0$, and the proof is complete. $\;\;\;\;\Box$

\begin{prop} There are finitely many irreducible
$A(L(n-l+\frac{1}{2},0))$-modules from the category $\mathcal{O}$.
\end{prop}
{\bf Proof:} It follows from Corollary \ref{c.1.7.2} that highest weights
$\mu \in {\frak h}^{*}$ of irreducible $A(L(n-l+\frac{1}{2},0))$-modules
$V(\mu)$ satisfy $p(\mu)=0$ for all $p \in {\mathcal P}_{0}$.
Lemma \ref{l.2.4.2} implies that $p_{1}(\mu)=p_{2}(\mu)=
\ldots =p_{l}(\mu)=0$ for such weights $\mu$.
Every weight $\mu \in {\frak h}^{*}$ is uniquely determined by its
values $\mu _{i}=\mu(h_{i})$, for $i=1, \ldots ,l$.
The equation $p_{l}(\mu)=0$ is
\begin{eqnarray*}
\mu _{l}(\mu _{l}-1)\cdot \ldots \cdot (\mu _{l}-2n+1)=0,
\end{eqnarray*}
which implies that there are $2n$
distinct values for $\mu _{l}$. The equation $p_{l-1}(\mu)=0$ is
\begin{eqnarray*}
 \mu_{l-1}\cdot \ldots \cdot (\mu_{l-1}-n+1)
(\mu_{l-1}+\mu_{l}+\frac{1}{2})\cdot \ldots
\cdot (\mu_{l-1}+\mu_{l}-n+\frac{3}{2})=0,
\end{eqnarray*}
which implies that for every fixed $\mu _{l}$, there are $2n$ distinct
values for $\mu _{l-1}$.
If we continue this procedure, at the end we get that equation
$p_{1}(\mu)=0$ is
\begin{eqnarray*}
&& \mu_{1}\cdot \ldots \cdot (\mu_{1}-n+1)
(\mu_{1}+ 2 \mu_{2}+ \ldots + 2 \mu_{l-1}+ \mu_{l}+l-\frac{3}{2})\cdot \ldots
\cdot \\
&& \quad \cdot (\mu_{1}+ 2 \mu_{2}+ \ldots + 2 \mu_{l-1}+
\mu_{l}+l-n-\frac{1}{2})=0,
\end{eqnarray*}
which implies that for fixed $\mu _{2}, \ldots ,\mu _{l}$, there are $2n$
distinct values for $\mu _{1}$.

Thus, there are at most $(2n)^{l}$ weights $\mu \in {\frak h}^{*}$
such that $p_{1}(\mu)=p_{2}(\mu)=
\ldots =p_{l}(\mu)=0$, which implies the claim of the theorem. $\;\;\;\;\Box$

It follows from Zhu's theory that:

\begin{thm} \label{c.2.4.3} There are finitely many irreducible weak
$L(n-l+\frac{1}{2},0)$-modules from the category $\mathcal{O}$.
\end{thm}

\subsection{Classification of irreducible $L(n-l+\frac{1}{2},0)$-modules}

In this section we classify irreducible
$L(n-l+\frac{1}{2},0)$-modules. It follows from Proposition \ref{t.1.3.5}
that irreducible ${\Z}_{+}$-graded
weak $L(n-l+\frac{1}{2},0)$-modules are in one-to-one correspondence with
irreducible $A(L(n-l+\frac{1}{2},0))$-modules. Specially,
irreducible $L(n-l+\frac{1}{2},0)$-modules are in
one-to-one correspondence with finite-dimensional
irreducible $A(L(n-l+\frac{1}{2},0))$-modules.

It follows from Proposition \ref{t.2.3.5} that
$A(L(n-l+\frac{1}{2},0)) \cong \frac{U(\frak g)}{I_{n}},$
which implies that every finite-dimensional
irreducible $A(L(n-l+\frac{1}{2},0))$-module is a
finite-dimensional irreducible $\frak g$-module, and
therefore is of the form $V(\mu)$, where $\mu \in P_{+}$ is a dominant
integral weight of $\frak g$.
\begin{lem} \label{t.2.4.4}
Assume that $V(\mu), \mu \in P_{+}$ is an
$A(L(n-l+\frac{1}{2},0))$-module. Then
$( \mu , \epsilon_{1} ) \leq n-\frac{1}{2}$.
\end{lem}
{\bf Proof:} If $V(\mu), \mu \in P_{+}$ is an
$A(L(n-l+\frac{1}{2},0))$-module, then Corollary
\ref{c.1.7.2} and Lemma \ref{l.2.4.2} imply that
$q(\mu)=0$. Thus,
\begin{eqnarray*}
& &0={\displaystyle \sum_{{(k_{1},\ldots ,k_{l}) \in {\Z}_{+}^{l} \atop \sum k_{i}=n}}} \frac{1}{k_{1}!4^{k_{1}}} \cdot
(\mu (h_{\epsilon_1})-2k_{2}-\ldots -2k_{l})\cdot \ldots \cdot
(\mu (h_{\epsilon_1})-2n+1)\cdot \\
& &(\mu (h_{\epsilon_1 - \epsilon_l})
-k_{l-1}-\ldots -k_{2}) \cdot \ldots \cdot (\mu (h_{\epsilon_1 -
\epsilon_l}) -k_{l-1}-\ldots -k_{2}-k_{l}-1) \\
& & \cdot \ldots \cdot
\mu (h_{\epsilon_1 - \epsilon_2}) \cdot \ldots \cdot (\mu (h_{\epsilon_1 -
\epsilon_2})-k_{2}+1).
\end{eqnarray*}

We claim that this relation implies
$( \mu , \epsilon_{1} ) \leq n-\frac{1}{2}$, or
equivalently $\mu (h_{\epsilon_{1}} ) \leq 2n-1$.
Suppose that $\mu (h_{\epsilon_{1}} ) \geq 2n$.
We claim that, under that assumption, all the summands above
are non-negative. Let $(k_{1},\ldots ,k_{l}) \in {\Z}_{+}^{l}$
be any $l$-tuple, such that $\sum_{i=1}^{l} k_{i}=n$. It is clear that
$$(\mu (h_{\epsilon_1})-2k_{2}-\ldots -2k_{l})\cdot \ldots \cdot
(\mu (h_{\epsilon_1})-2n+1) \geq 0.$$
Assume that
\begin{eqnarray*}
& &(\mu (h_{\epsilon_1 - \epsilon_l})
-k_{l-1}-\ldots -k_{2}) \cdot \ldots \cdot (\mu (h_{\epsilon_1 -
\epsilon_l}) -k_{l-1}-\ldots -k_{2}-k_{l}-1) \\
& & \cdot \ldots \cdot
\mu (h_{\epsilon_1 - \epsilon_2}) \cdot \ldots \cdot (\mu (h_{\epsilon_1 -
\epsilon_2})-k_{2}+1) \neq 0.
\end{eqnarray*}
Then, from
$\mu (h_{\epsilon_1 - \epsilon_2}) \cdot \ldots \cdot (\mu (h_{\epsilon_1 -
\epsilon_2})-k_{2}+1)\neq 0,$
and from $\mu (h_{\epsilon_1 - \epsilon_2}) \in {\Z}_{+}$
follows $\mu (h_{\epsilon_1 - \epsilon_2}) \geq k_{2},$
which implies
$\mu (h_{\epsilon_1 - \epsilon_2}) \cdot \ldots \cdot (\mu (h_{\epsilon_1 -
\epsilon_2})-k_{2}+1) > 0.$
Since $h_{\epsilon_1 - \epsilon_3}=h_{\epsilon_1 - \epsilon_2}
+h_{\epsilon_2 - \epsilon_3}$, we have
$\mu (h_{\epsilon_1 - \epsilon_3}) \geq
\mu (h_{\epsilon_1 - \epsilon_2}) \geq k_{2}.$
From
$\mu (h_{\epsilon_1 - \epsilon_3}-k_{2}) \cdot \ldots \cdot
(\mu (h_{\epsilon_1 -\epsilon_3})-k_{2}-k_{3}+1)\neq 0,$
and from $\mu (h_{\epsilon_1 - \epsilon_3}) \geq k_{2}$
follows $\mu (h_{\epsilon_1 - \epsilon_3}) \geq k_{2}+k_{3},$
which implies
$\mu (h_{\epsilon_1 - \epsilon_3}-k_{2}) \cdot \ldots \cdot
(\mu (h_{\epsilon_1 -\epsilon_3})-k_{2}-k_{3}+1)>0.$
Inductively, we obtain
\begin{eqnarray*}
& &(\mu (h_{\epsilon_1 - \epsilon_l})
-k_{l-1}-\ldots -k_{2}) \cdot \ldots \cdot (\mu (h_{\epsilon_1 -
\epsilon_l}) -k_{l-1}-\ldots -k_{2}-k_{l}-1) \\
& & \cdot \ldots \cdot
\mu (h_{\epsilon_1 - \epsilon_2}) \cdot \ldots \cdot (\mu (h_{\epsilon_1 -
\epsilon_2})-k_{2}+1) > 0.
\end{eqnarray*}
Thus, we have shown that all the summands above are non-negative.

On the other hand, for $k_{1}=n$, $k_{2}=k_{3}= \ldots =k_{l}=0$,
we have the summand
$$\frac{1}{n!4^{n}}
\mu (h_{\epsilon_1}) \cdot \ldots \cdot
(\mu (h_{\epsilon_1})-2n+1)>0.$$
Therefore
\begin{eqnarray*}
& &0<{\displaystyle \sum_{{(k_{1},\ldots ,k_{l}) \in {\Z}_{+}^{l} \atop \sum k_{i}=n}}} \frac{1}{k_{1}!4^{k_{1}}} \cdot
(\mu (h_{\epsilon_1})-2k_{2}-\ldots -2k_{l})\cdot \ldots \cdot
(\mu (h_{\epsilon_1})-2n+1)\cdot \\
& &(\mu (h_{\epsilon_1 - \epsilon_l})
-k_{l-1}-\ldots -k_{2}) \cdot \ldots \cdot (\mu (h_{\epsilon_1 -
\epsilon_l}) -k_{l-1}-\ldots -k_{2}-k_{l}-1) \\
& & \cdot \ldots \cdot
\mu (h_{\epsilon_1 - \epsilon_2}) \cdot \ldots \cdot (\mu (h_{\epsilon_1 -
\epsilon_2})-k_{2}+1),
\end{eqnarray*}
which is a contradiction. Thus, $( \mu , \epsilon_{1} ) \leq
n-\frac{1}{2}$. $\;\;\;\;\Box$

The converse of Lemma \ref{t.2.4.4} also holds:

\begin{lem} Let $\mu \in P_{+}$, such that
$( \mu , \epsilon_{1} ) \leq n-\frac{1}{2}$. Then
$V(\mu)$ is an $A(L(n-l+\frac{1}{2},0))$-module.
\end{lem}
{\bf Proof:} Since
$A(L(n-l+\frac{1}{2},0)) \cong \frac{U(\frak g)}{I_{n}},$
we have to show $I_{n}.V(\mu)=0$. Since the ideal $I_{n}$ is generated
by the vector $v_n'$, it is sufficient to show that $v_n'$
annihilates $V(\mu)$.

Suppose that there exists a vector $u \in V(\mu)$,
such that $v_n'.u \neq 0$. From the structure of the root
system of type $B_{l}$ follows that the lowest weight in the
module $V(\mu)$ is $- \mu$.
Then the weight of the vector $u$ is of the form
$- \mu + \sum_{i=1}^{l} k_{i} \alpha_{i}$, where
$k_{i} \in {\Z}_{+}$ for $i=1,\ldots ,l$,
and the weight of the vector $v_n'.u$ has to be
of the form $ \mu - \sum_{i=1}^{l} t_{i} \alpha_{i}$,
where $t_{i} \in {\Z}_{+}$ for $i=1,\ldots ,l$.
Since $v_n'$ has the weight
$2n \epsilon_{1}$, we obtain the equation:
$- \mu + \sum_{i=1}^{l} k_{i} \alpha_{i} + 2n \epsilon_{1}=
\mu - \sum_{i=1}^{l} t_{i} \alpha_{i},$
or equivalently
$$2 \mu - 2n \epsilon_{1}= \sum_{i=1}^{l} m_{i} \alpha_{i},$$
where $m_{i}=k_{i}+t_{i} \in {\Z}_{+}$ for $i=1,\ldots ,l$.
It follows from the equation above that
$(2 \mu - 2n \epsilon_{1}, \epsilon_{1})=
(\sum_{i=1}^{l} m_{i} \alpha_{i}, \epsilon_{1}).$
Since $( \alpha_{1}, \epsilon_{1})=1$ and
$( \alpha_{i}, \epsilon_{1})=0$ for $i=2,\ldots ,l$,
we obtain
$2 (\mu , \epsilon_{1})-2n=m_{1}.$
From $( \mu , \epsilon_{1} ) \leq n-\frac{1}{2}$
we obtain
$$m_{1}=2 (\mu , \epsilon_{1})-2n \leq -1,$$
which is a contradiction with $m_{i} \in {\Z}_{+}$
for $i=1,\ldots ,l$.

Thus, $v_n'$ annihilates $V(\mu)$, and
$V(\mu)$ is an $A(L(n-l+\frac{1}{2},0))$-module. $\;\;\;\;\Box$

\begin{prop} The set
$$ \{ V(\mu) \ \vert \ \mu \in P_{+}, \ ( \mu , \epsilon_{1} ) \leq
n-\frac{1}{2} \}$$
provides the complete list of irreducible finite-dimensional
$A(L(n-l+\frac{1}{2},0))$-modules.
\end{prop}

It follows from Zhu's theory that:

\begin{thm} \label{k.2.4.6}
The set
$$ \{ L(n-l+\frac{1}{2}, \mu) \ \vert \ \mu \in P_{+}, \ ( \mu , \epsilon_{1} ) \leq
n-\frac{1}{2} \}$$
provides the complete list of irreducible
$L(n-l+\frac{1}{2},0)$-modules.
\end{thm}

\subsection{Complete reducibility in category of $L(n-l+\frac{1}{2},0)$-modules}

In this subsection we show that every $L(n-l+\frac{1}{2},0)$-module
is completely reducible. The following lemma is
crucial for proving complete reducibility.

\begin{lem} \label{p.2.5.1}
Let $L(\lambda)$ be a $L(n-l+\frac{1}{2},0)$-module.
Then the weight $\lambda$ is admissible.
\end{lem}
{\bf Proof:} If $L(\lambda)$ is a $L(n-l+\frac{1}{2},0)$-module,
then Theorem \ref{k.2.4.6} implies that $\lambda=(n-l+\frac{1}{2})
\Lambda_{0}+ \mu$, for some weight $\mu \in P_{+}$,
such that  $( \mu , \epsilon_{1} ) \leq
n-\frac{1}{2}$. It follows that
\begin{eqnarray*}
&& \langle \lambda + \rho,\alpha _{i}^{\vee}\rangle =
\langle \mu ,\alpha _{i}^{\vee}\rangle +1 \in \N \ \ \mbox{for } i=1,
\ldots ,l, \\
&& \langle \lambda + \rho,(\delta - \epsilon_{1})^{\vee} \rangle =
2n-2( \mu , \epsilon_{1} ) \in \N,
\end{eqnarray*}
which implies that
$\lambda$ is admissible weight and that
$ \hat{\Pi}^{\vee}_{\lambda }= \{ (\delta - \epsilon_{1})^{\vee},
\alpha_{1}^{\vee},\alpha_{2}^{\vee}, \ldots ,
\alpha_{l}^{\vee} \}.$ $\Box$

\begin{lem} \label{p.2.5.1.1}
Let $M$ be a $L(n-l+\frac{1}{2},0)$-module.
Then $M$ is from the category $\mathcal{O}$ as a $\hat{\frak g}$-module.
\end{lem}
{\bf Proof:} Let $M=\oplus_{\alpha\in {\C}}M_{(\alpha)}$,
where $L(0)u=\alpha u$ for any $u\in M_{(\alpha)}$,
$\dim M_{(\alpha)}<\infty$ for any $\alpha\in
{\C}$ and $M_{(\alpha+n)}=0$ for $n\in {\Z}$ sufficiently
small. It follows from
$a_mM_{(\alpha)}\subset M_{(\alpha+ {\rm wt} a-m-1)}$ for
$a \in V$, that $M_{(\alpha)}$ is a ${\frak g}$-module, for any
$\alpha\in {\C}$. Since $M_{(\alpha)}$ is
finite-dimensional, $\frak h$ acts semisimply
on $M_{(\alpha)}$, which implies that $\hat{\frak h}$
acts semisimply on $M$ with finite-dimensional weight spaces. Let $v \in M$ be a
singular vector of weight $\lambda \in \hat{\frak h} ^* $.
Then $L(\lambda)$ is an irreducible subquotient
of $M$, which implies that $L(\lambda)$ is a $L(n-l+\frac{1}{2},0)$-module.
It follows from Theorem \ref{k.2.4.6} that there are finitely many
irreducible \linebreak $L(n-l+\frac{1}{2},0)$-modules, which implies that
there exists a finite number of weights $\nu_{1}, \ldots , \nu_{k}
\in \hat{{\frak h}}^{*}$ such that for every weight
$\nu$ of $M$ holds $\nu \in \cup_{i=1}^{k}D(\nu_{i})$. Thus
$\hat{\frak g}$-module $M$ is from the category $\mathcal{O}$.
$\;\;\;\;\Box$

\begin{thm} \label{t.2.5.2}
Let $M$ be a $L(n-l+\frac{1}{2},0)$-module. Then
$M$ is completely reducible.
\end{thm}
{\bf Proof:} Let $L(\lambda)$ be some irreducible subquotient
of $M$. Then $L(\lambda)$ is a $L(n-l+\frac{1}{2},0)$-module,
and Lemma \ref{p.2.5.1} implies that $\lambda$ is admissible
weight. It follows from Lemma \ref{p.2.5.1.1} that $M$ is
from the category $\mathcal{O}$, and then Proposition \ref{t.KW2} implies
that $M$ is completely reducible. $\;\;\;\;\Box$

\section{Weak $L(-l+\frac{3}{2},0)$-modules
from category~$\mathcal{O}$}

In this section we study the special case $n=1$, i.e.
the smallest admissible half-integer level $-l+\frac{3}{2}$.
In this case we find a basis for the vector space
$ {\mathcal P}_{0}$, defined in subsection \ref{subsec.3.3},
from which we get the classification of irreducible weak
$L(-l+\frac{3}{2},0)$-modules from the category $\mathcal{O}$.
We also show that every weak
$L(-l+\frac{3}{2},0)$-module from the category $\mathcal{O}$
is completely reducible.

\subsection{Classification of irreducible weak
$L(-l+\frac{3}{2},0)$-modules from category~$\mathcal{O}$}

It follows from Corollary \ref{c.1.7.2} that irreducible
$A(L(-l+\frac{3}{2},0))$-modules are in one-to-one
correspondence with weights $\mu \in {\frak h}^{*}$ such that
$p(\mu)=0$ for all $p \in {\mathcal P}_{0}$, where
$ {\mathcal P}_{0}=\{ \ p_{r} \ \vert \ r \in R_{0} \}.$
In this case, $R$ is a highest weight $U(\frak g)$-module with the
highest weight $2 \epsilon_{1}=2 \omega_{1}$, and $R_{0}$ is zero-weight
subspace of $R$.

\begin{lem} \label{l.2.6.0}
$$ \dim R_{0} \leq l $$
\end{lem}
{\bf Proof:} In this proof we use induction on $l$. We use the
notation $V_{l}(\mu)$ for the highest weight module
for simple Lie algebra of type $B_{l}$, with the highest
weight $\mu \in {\frak h}^{*}$.

For $l=2$ it is easily checked that $ \dim R =14$ and $ \dim R_{0} =2$.
Assume that the claim of this lemma holds for simple Lie algebra
of type $B_{l-1}$, $l-1 \geq 2$. Let $\frak g$ be simple Lie algebra
of type $B_{l}$.
Let $\frak g'$ be the subalgebra of $\frak g$ associated to roots
$\alpha_{2}, \ldots ,\alpha_{l}$. $\frak g'$ is then a simple Lie algebra
of type $B_{l-1}$. We can decompose $\frak g$-module $V_{l}(2 \omega_{1})$
into the direct sum of irreducible $\frak g'$-modules. If we denote by $v$ the
highest weight vector of $\frak g$-module $V_{l}(2 \omega_{1})$,
then it can be easily checked that $f_{\epsilon_1 - \epsilon_2}^{2}.v$,
$f_{\epsilon_1 - \epsilon_2}.v$, $f_{\epsilon_1}^{2}
f_{\epsilon_1 - \epsilon_2}.v$, $v$, $f_{\epsilon_1}^{4}.v$ and
$(-\frac{1}{4}f_{\epsilon_1}^2+ f_{\epsilon_1 - \epsilon_2}
f_{\epsilon_1 + \epsilon_2}+ \cdots +f_{\epsilon_1 - \epsilon_l}
f_{\epsilon_1 + \epsilon_l}).v$ are highest weight vectors for
$\frak g'$, which generate $\frak g'$-modules isomorphic to
$V_{l-1}(2 \omega_{1})$, $V_{l-1}(\omega_{1})$,
$V_{l-1}(\omega_{1})$, $V_{l-1}(0)$, $V_{l-1}(0)$ and $V_{l-1}(0)$,
respectively. It follows from Weyl's formula for the dimension
of irreducible module that $ \dim V_{l}(2 \omega_{1})= \linebreak 2l^{2}+3l$
and $ \dim V_{l}(\omega_{1})=2l+1$, which implies that the direct
sum of $\frak g'$-modules above is $V_{l}(2 \omega_{1})$.
Clearly, there are no zero-weight vectors for $\frak g$ in
$\frak g'$-modules generated by highest weight vectors
$f_{\epsilon_1 - \epsilon_2}.v$, $f_{\epsilon_1}^{2}
f_{\epsilon_1 - \epsilon_2}.v$, $v$ and $f_{\epsilon_1}^{4}.v$.
The inductive assumption implies that there are at most $l-1$
linearly independent zero-weight vectors for $\frak g$ in
$\frak g'$-module $V_{l-1}(2 \omega_{1})$, which implies
$ \dim R_{0} \leq (l-1)+1=l$. $\;\;\;\;\Box$

In the next lemma we find a basis for the vector space
$ {\mathcal P}_{0}$.

\begin{lem} \label{l.2.6.1}
\begin{eqnarray*}
{\mathcal P}_{0}= \mbox{span}_{\C} \{ p_{1}, \ldots ,p_{l} \},
\end{eqnarray*}
where
\begin{eqnarray*}
&& p_{i}(h)=h_{i}(h_{\epsilon_i + \epsilon_{i+1}}+l-i-\frac{1}{2}), \
\mbox{for } i=1, \ldots ,l-1, \\
&& p_{l}(h)=h_{l}(h_{l}-1).
\end{eqnarray*}
\end{lem}
{\bf Proof:} Lemma \ref{l.2.4.2} implies that
$p_{1}, \ldots ,p_{l}$ are linearly independent polynomials
in the set ${\mathcal P}_{0}$. It follows from the definition
of the set ${\mathcal P}_{0}$ that $\dim {\mathcal P}_{0} \leq \dim R_{0}$,
and using Lemma \ref{l.2.6.0} we get $\dim {\mathcal P}_{0} \leq l$.
Thus, polynomials $p_{1}, \ldots ,p_{l}$ form a basis for
${\mathcal P}_{0}$. $\;\;\;\;\Box$

\begin{prop}\label{t.2.6.3}
For every subset $S=\{ i_{1}, \ldots ,i_{k} \}
\subseteq \{1,2, \ldots , l-1 \}$, \linebreak $i_{1}< \ldots
<i_{k}$, we define weights
\begin{eqnarray*}
&& \mu _{S}= \sum _{j=1}^{k}\left(i_{j}+2 \sum _{s=j+1}^{k}
(-1)^{s-j}i_{s} + (-1)^{k-j+1}(l- \frac{1}{2})\right) \omega _{i_{j}}, \\
&& \mu _{S}'=\sum _{j=1}^{k}\left(i_{j}+2 \sum _{s=j+1}^{k}
(-1)^{s-j}i_{s} + (-1)^{k-j+1}(l+ \frac{1}{2})\right) \omega _{i_{j}}
+\omega _{l},
\end{eqnarray*}
where $\omega _{1}, \ldots , \omega _{l}$ are fundamental
weights for $\frak g$. Then the set
$$ \{ V(\mu _{S}), V(\mu _{S}') \ \vert \ S \subseteq \{1,2, \ldots , l-1
\} \}$$
provides the complete list of irreducible $A(L(-l+\frac{3}{2},0))$-modules
from the category $\mathcal{O}$.
\end{prop}
{\bf Proof:} Since
$$ h_{\epsilon_i + \epsilon_{i+1}}=h_{i}+2h_{i+1}+ \ldots
+2h_{l-1}+h_{l}, $$
it follows from Corollary \ref{c.1.7.2} and Lemma \ref{l.2.6.1}
that highest weights $\mu \in {\frak h}^{*}$ of irreducible
$A(L(-l+\frac{3}{2},0))$-modules $V(\mu)$ are in one-to-one
correspondence with solutions of the system of polynomial
equations
\begin{eqnarray*}
&& h_{1}(h_{1}+2h_{2}+ \ldots +2h_{l-1}+h_{l}+l-\frac{3}{2})=0,  \\
&& h_{2}(h_{2}+2h_{3}+ \ldots +2h_{l-1}+h_{l}+l-\frac{5}{2})=0,  \\
&& \qquad \qquad \qquad \vdots \\
&& h_{l-1}(h_{l-1}+h_{l}+\frac{1}{2})=0,  \\
&& h_{l}(h_{l}-1)=0.
\end{eqnarray*}
Clearly, $h_{l} \in \{ 0,1 \}$. Let
$S=\{ i_{1}, \ldots ,i_{k} \}$, $i_{1}< \ldots <i_{k}$ be the
subset of $\{1,2, \ldots , l-1 \}$ such that
$h_{i}=0$ for $i \notin S$ and $h_{i} \neq 0$ for $i \in S$.

First consider the case $h_{l}=0$. Then we have the system
\begin{eqnarray}
&& h_{i_{1}}+2h_{i_{2}}+ \ldots +2h_{i_{k}}+l-i_{1}-\frac{1}{2}=0, \nno \\
&& h_{i_{2}}+2h_{i_{3}}+ \ldots +2h_{i_{k}}+l-i_{2}-\frac{1}{2}=0, \nno \\
&& \qquad \qquad \qquad \vdots \label{2.6.1.1}\\
&& h_{i_{k-1}}+2h_{i_{k}}+l-i_{k-1}-\frac{1}{2}=0, \nno \\
&& h_{i_{k}}+l-i_{k}-\frac{1}{2}=0. \nno
\end{eqnarray}
The solution of this triangular system is
$$h_{i_{j}}= i_{j}+2 \sum _{s=j+1}^{k}
(-1)^{s-j}i_{s} + (-1)^{k-j+1}(l- \frac{1}{2}),
\ \mbox{ for } j=1, \ldots ,k.$$
It follows that $V(\mu _{S})$ is irreducible
$A(L(-l+\frac{3}{2},0))$-module.

In the case when $h_{l}=1$, we have the system
\begin{eqnarray*}
&& h_{i_{1}}+2h_{i_{2}}+ \ldots +2h_{i_{k}}+l-i_{1}+\frac{1}{2}=0,  \\
&& h_{i_{2}}+2h_{i_{3}}+ \ldots +2h_{i_{k}}+l-i_{2}+\frac{1}{2}=0,  \\
&& \qquad \qquad \qquad \vdots \\
&& h_{i_{k-1}}+2h_{i_{k}}+l-i_{k-1}+\frac{1}{2}=0, \\
&& h_{i_{k}}+l-i_{k}+\frac{1}{2}=0,
\end{eqnarray*}
whose solution is
$$h_{i_{j}}= i_{j}+2 \sum _{s=j+1}^{k}
(-1)^{s-j}i_{s} + (-1)^{k-j+1}(l+ \frac{1}{2}),
\ \mbox{ for } j=1, \ldots ,k.$$
It follows that $V(\mu _{S}')$ is irreducible
$A(L(-l+\frac{3}{2},0))$-module,
which completes the proof. $\;\;\;\;\Box$

It follows from Zhu's theory that:

\begin{thm} \label{c.2.6.3}
The set
$$ \{ L(-l+\frac{3}{2}, \mu _{S}), L(-l+\frac{3}{2}, \mu _{S}') \ \vert \
S \subseteq \{1,2, \ldots , l-1 \} \}$$
provides the complete list of irreducible weak
$L(-l+\frac{3}{2},0)$-modules from the category $\mathcal{O}$.
\end{thm}

Theorem \ref{c.2.6.3} implies that there are $2^{l}$ irreducible weak
$L(-l+\frac{3}{2},0)$-modules from the category $\mathcal{O}$.

\subsection{Complete reducibility of weak $L(-l+\frac{3}{2},0)$-modules
from category $\mathcal{O}$} \label{sec.2.7}

We introduce the notation $\lambda _{S}=(-l+\frac{3}{2}) \Lambda_{0}+ \mu _{S}$
and $\lambda _{S}'=(-l+\frac{3}{2}) \Lambda_{0}+ \mu _{S}'$, for every
$S \subseteq \{1,2, \ldots , l-1 \}$. The following lemma is
crucial for proving complete reducibility.

\begin{lem}\label{l.2.6.4}
Weights  $\lambda _{S}, \lambda _{S}' \in \hat{\frak h} ^{*}$
are admissible, for every $S \subseteq \{1,2, \ldots , l~-~1 \}$.
\end{lem}
{\bf Proof:} We will prove that the weight $\lambda _{S}$
is admissible for every $S \subseteq \{1,2, \ldots , l-1 \}$. The
proof of admissibility of weights $\lambda _{S}'$
for every $S \subseteq \{1,2, \ldots , l-1 \}$ is similar.
We have to show
\begin{eqnarray}
& &\langle \lambda_{S} + \rho,\tilde{\alpha}^{\vee}\rangle \notin -{\Z}_{+} \mbox{ for
any }
\tilde{\alpha} \in \hat{\Delta}^{\mbox{\scriptsize{re}}}_{+}, \label{2.6.3.1}\\
& &{\Q} \hat{\Delta}^{\vee \mbox{\scriptsize{re}}}_{\lambda _{S}}={\Q}
\hat{\Pi}^{\vee} .\label{2.6.3.2}
\end{eqnarray}
Let's prove first the relation (\ref{2.6.3.1}). Any positive real root
$\tilde{\alpha} \in \hat{\Delta}^{\mbox{\scriptsize{re}}}_{+}$
of $\hat{\frak g}$ is of the form $\tilde{\alpha}= \alpha +m \delta$, for $m>0$ and $\alpha
\in \Delta$ or $m=0$ and $\alpha \in {\Delta}_{+}$. It can be easily
checked, as in Lemma \ref{l.2.3.1}, that
\begin{eqnarray}
& & \langle \lambda_{S} + \rho,\tilde{\alpha}^{\vee}\rangle =
\frac{2}{( \alpha , \alpha )}\left(m(l+\frac{1}{2})+( \bar{\rho}
, \alpha )+ ( \mu _{S} , \alpha ) \right), \label{2.6.2.1}
\end{eqnarray}
where $\bar{\rho}$ is the sum of fundamental weights of $\frak g$.
Let
$S=\{ i_{1}, \ldots ,i_{k} \} \subseteq \{1,2, \ldots , l-1 \}$,
$i_{1}< \ldots <i_{k}$. Proposition \ref{t.2.6.3} implies that
$ \mu _{S}= \sum _{j=1}^{k} h_{i_{j}} \omega _{i_{j}}$,
where
$$h_{i_{j}}= i_{j}+2 \sum _{s=j+1}^{k}
(-1)^{s-j}i_{s} + (-1)^{k-j+1}(l- \frac{1}{2}),
\ \mbox{ for } j=1, \ldots ,k.$$
From the system (\ref{2.6.1.1}) easily follows
\begin{eqnarray}
h_{i_{j}}+h_{i_{j+1}}=-(i_{j+1}-i_{j}) \label{2.6.2}\
\mbox{ for } j=1, \ldots ,k-1.
\end{eqnarray}
We will consider three cases in this proof. \\
{\bf Case 1.}: The root $\alpha$ is of the form
$\alpha= \pm  (\epsilon_{i} - \epsilon_{j})$, for $i,j=1, \ldots ,l$,
$i<j$.

Then $( \bar{\rho}, \epsilon_{i} -
\epsilon_{j})=j-i$. Let $s,t \in \{ 1, \ldots ,k \}$ be indices
such that $S \cap \{ i,i+1, \ldots ,j-1 \}= \{ i_{s}, \ldots ,i_{t}
\}$. Clearly, $i_{s} \geq i$ and $i_{t} \leq j-1$. Furthermore,
$( \mu _{S}, \epsilon_{i} -
\epsilon_{j})=h_{i_{s}}+ \ldots \ +h_{i_{t}}$.

First consider the case
$\alpha= \epsilon_{i} - \epsilon_{j}$, $i<j$ and $m \geq 0$.

If $t-s+1$ is even, then using relation
(\ref{2.6.2}) we obtain
\begin{eqnarray*}
( \mu _{S}, \epsilon_{i} - \epsilon_{j})=
-(i_{s+1}-i_{s})- \ldots - (i_{t}-i_{t-1}) \geq -(i_{t}-i_{s}),
\end{eqnarray*}
and relation (\ref{2.6.2.1}) implies
\begin{eqnarray*}
&& \langle \lambda_{S} + \rho,\tilde{\alpha}^{\vee}\rangle \geq
( \bar{\rho} , \epsilon_{i} - \epsilon_{j} )+ ( \mu _{S} , \epsilon_{i} - \epsilon_{j} )
\geq (j-i)- (i_{t}-i_{s})= \\
&& \qquad \qquad \quad \ \ =(j-i_{t})+ (i_{s}-i) >0.
\end{eqnarray*}

Suppose now that $t-s+1$ is odd. Then
$( \mu _{S}, \epsilon_{i} - \epsilon_{j}) \notin \Z$,
and if $m=0$, then $\langle \lambda_{S} + \rho,\tilde{\alpha}^{\vee}\rangle
\notin \Z$. Let $m \geq 1$. Then
\begin{eqnarray*}
&& ( \mu _{S}, \epsilon_{i} - \epsilon_{j})=h_{i_{s}}+ \ldots \ +h_{i_{t-1}}
+h_{i_{t}}= -(i_{s+1}-i_{s})- \ldots -
(i_{t-1}-i_{t-2})+h_{i_{t}} \\
&& \qquad \qquad \quad \  \geq -(i_{t-1}-i_{s})+h_{i_{t}}.
\end{eqnarray*}
If $k-t$ is even, we get
\begin{eqnarray*}
h_{i_{t}}=i_{t}+ 2(i_{t+2}-i_{t+1})+ \ldots + 2(i_{k}-i_{k-1})
-(l- \frac{1}{2}) \geq -(l- \frac{1}{2}),
\end{eqnarray*}
which implies
\begin{eqnarray*}
( \mu _{S}, \epsilon_{i} - \epsilon_{j}) \geq -(i_{t-1}-i_{s})-(l- \frac{1}{2}).
\end{eqnarray*}
We obtain
\begin{eqnarray*}
&& \langle \lambda_{S} + \rho,\tilde{\alpha}^{\vee}\rangle \geq
l+ \frac{1}{2} +( \bar{\rho} , \epsilon_{i} - \epsilon_{j} )+ ( \mu _{S} ,
\epsilon_{i} - \epsilon_{j}) \\
&& \qquad \qquad \quad \ \geq l+ \frac{1}{2}+ (j-i)-(i_{t-1}-i_{s})-(l- \frac{1}{2}) \\
&& \qquad \qquad \quad \ = (j-i_{t-1})+ (i_{s}-i)+1 >0.
\end{eqnarray*}
If $k-t$ is odd, then
\begin{eqnarray*}
h_{i_{t}}=i_{t}+ 2(i_{t+2}-i_{t+1})+ \ldots + 2(i_{k-1}-i_{k-2})-2i_{k}
+(l- \frac{1}{2}) \geq l- \frac{1}{2} -2i_{k},
\end{eqnarray*}
which implies
\begin{eqnarray*}
( \mu _{S}, \epsilon_{i} - \epsilon_{j}) \geq -(i_{t-1}-i_{s})+ l- \frac{1}{2} -2i_{k}.
\end{eqnarray*}
We obtain
\begin{eqnarray*}
&& \langle \lambda_{S} + \rho,\tilde{\alpha}^{\vee}\rangle \geq
l+ \frac{1}{2} +( \bar{\rho} , \epsilon_{i} - \epsilon_{j} )+ ( \mu _{S} , \epsilon_{i} - \epsilon_{j} )
\geq l+ \frac{1}{2}+ (j-i)- \\
&& \quad  \quad  -(i_{t-1}-i_{s})+l- \frac{1}{2} -2i_{k} = 2(l-i_{k})+(j-i_{t-1})+ (i_{s}-i) >0.
\end{eqnarray*}

Thus, we have proved that, if
$\alpha= \epsilon_{i} - \epsilon_{j}$, $i<j$ and $m \geq 0$,
then \linebreak $ \langle \lambda_{S} + \rho,\tilde{\alpha}^{\vee}\rangle \notin
-{\Z}_{+}$.

Now, let's consider the case
$\alpha= -(\epsilon_{i} - \epsilon_{j})$, $i<j$ and $m \geq 1$.

Then
\begin{eqnarray*}
& & \langle \lambda_{S} + \rho,\tilde{\alpha}^{\vee}\rangle =
m(l+\frac{1}{2})-( \bar{\rho}, \epsilon_{i} - \epsilon_{j} )-
( \mu _{S} , \epsilon_{i} - \epsilon_{j} ).
\end{eqnarray*}

If $t-s+1$ is even, then $( \mu _{S}, \epsilon_{i} -
\epsilon_{j})$ is an integer and $( \mu _{S}, \epsilon_{i} -
\epsilon_{j}) \leq 0$, so if $m$ is odd,
then $\langle \lambda_{S} + \rho,\tilde{\alpha}^{\vee}\rangle
\notin \Z$. If $m$ is even, then $m \geq 2$, and we get
\begin{eqnarray*}
\langle \lambda_{S} + \rho,\tilde{\alpha}^{\vee}\rangle \geq
2(l+\frac{1}{2})-(j-i)=(l-j)+l+i+1>0.
\end{eqnarray*}
If $t-s+1$ is odd, then
$( \mu _{S}, \epsilon_{i} - \epsilon_{j})=h_{i_{s}}+ \ldots \ +h_{i_{t-1}}
+h_{i_{t}} \leq h_{i_{t}}$.
If $k-t$ is even, then
\begin{eqnarray*}
&& h_{i_{t}}=i_{t}+ 2(i_{t+2}-i_{t+1})+ \ldots + 2(i_{k}-i_{k-1})
-(l- \frac{1}{2}) \leq  \\
&&  \quad \ \ \leq i_{t}+ 2(i_{k}-i_{t+1})-(l-
\frac{1}{2}),
\end{eqnarray*}
which implies
\begin{eqnarray*}
( \mu _{S}, \epsilon_{i} - \epsilon_{j}) \leq i_{t}+ 2(i_{k}-i_{t+1})-(l-
\frac{1}{2}).
\end{eqnarray*}
It follows
\begin{eqnarray*}
&& \langle \lambda_{S} + \rho,\tilde{\alpha}^{\vee}\rangle \geq
l+ \frac{1}{2} -(j-i)- \left( i_{t}+ 2(i_{k}-i_{t+1})-(l-
\frac{1}{2}) \right) \\
&& \qquad \qquad \quad \ =
2(l-i_{k})+(i_{t+1}-i_{t})+(i_{t+1}-j)+i>0.
\end{eqnarray*}
If $k-t$ is odd, then
\begin{eqnarray*}
&& h_{i_{t}}=i_{t}+ 2(i_{t+2}-i_{t+1})+ \ldots + 2(i_{k-1}-i_{k-2})
-2i_{k}+(l- \frac{1}{2}) \leq  \\
&&  \quad \ \ \leq i_{t}+ 2(i_{k-1}-i_{t+1})-2i_{k}+(l- \frac{1}{2}),
\end{eqnarray*}
which implies
\begin{eqnarray*}
( \mu _{S}, \epsilon_{i} - \epsilon_{j}) \leq i_{t}+ 2(i_{k-1}-i_{t+1})-2i_{k}+(l- \frac{1}{2}).
\end{eqnarray*}
It follows
\begin{eqnarray*}
&& \langle \lambda_{S} + \rho,\tilde{\alpha}^{\vee}\rangle \geq
l+ \frac{1}{2} -(j-i)- \left( i_{t}+ 2(i_{k-1}-i_{t+1})-2i_{k}+(l- \frac{1}{2}) \right) \\
&& \qquad \qquad \quad \ =
2(i_{k}-i_{k-1})+(i_{t+1}-i_{t})+(i_{t+1}-(j-1))+i>0.
\end{eqnarray*}
We have proved that, if
$\alpha=- (\epsilon_{i} - \epsilon_{j})$, $i<j$ and $m \geq 1$,
then \linebreak $ \langle \lambda_{S} + \rho,\tilde{\alpha}^{\vee}\rangle \notin
-{\Z}_{+}$. \\
{\bf Case 2.}: The root $\alpha$ is of the form
$\alpha= \pm  \epsilon_{i}$, for $i=1, \ldots ,l$.

Then $( \bar{\rho}, \epsilon_{i})=l-i+\frac{1}{2}$. Let $s \in \{ 1, \ldots ,k \}$
be the index such that $S \cap \{ i,i+1, \ldots ,l-1 \}= \{ i_{s}, \ldots ,i_{k}
\}$. Clearly, $i_{s-1} < i \leq i_{s}$. Furthermore,
$( \mu _{S}, \epsilon_{i} )=h_{i_{s}}+ \ldots \ +h_{i_{k}}$.

First consider the case
$\alpha= \epsilon_{i}$ and $m \geq 0$.

Then
\begin{eqnarray}
& & \langle \lambda_{S} + \rho,\tilde{\alpha}^{\vee}\rangle =
2\left(m(l+\frac{1}{2})+( \bar{\rho}
, \epsilon_{i} )+ ( \mu _{S} , \epsilon_{i} ) \right).
\end{eqnarray}

If $k-s+1$ is even, then using relation
(\ref{2.6.2}) we obtain
\begin{eqnarray*}
( \mu _{S}, \epsilon_{i})=
-(i_{s+1}-i_{s})- \ldots - (i_{k}-i_{k-1}) \geq -(i_{k}-i_{s}),
\end{eqnarray*}
which implies
\begin{eqnarray*}
&& \langle \lambda_{S} + \rho,\tilde{\alpha}^{\vee}\rangle \geq
2( \bar{\rho} , \epsilon_{i} )+ 2( \mu _{S} , \epsilon_{i} )
\geq (2l-2i+1)- 2(i_{k}-i_{s}) \\
&& \qquad \qquad \quad \ = 2(l-i_{k})+ 2(i_{s}-i)+1 >0.
\end{eqnarray*}

If $k-s+1$ is odd, then
\begin{eqnarray*}
&& ( \mu _{S}, \epsilon_{i})=h_{i_{s}}+ \ldots \ +h_{i_{k-1}}
+h_{i_{k}}= -(i_{s+1}-i_{s})- \ldots -
(i_{k-1}-i_{k-2})+h_{i_{k}} \\
&& \qquad \qquad \quad \  \geq -(i_{k-1}-i_{s})+(-l+i_{k}+\frac{1}{2})
= -l+\frac{1}{2}+i_{k}-i_{k-1}+i_{s},
\end{eqnarray*}
which implies
\begin{eqnarray*}
&& \langle \lambda_{S} + \rho,\tilde{\alpha}^{\vee}\rangle \geq
2( \bar{\rho} , \epsilon_{i} )+ 2( \mu _{S} , \epsilon_{i} )
\geq (2l-2i+1)-l+\frac{1}{2}+i_{k}-i_{k-1}+i_{s} \\
&& \qquad \qquad \quad \ = (l-i)+ (i_{k}-i_{k-1})+(i_{s}-i)+\frac{3}{2} >0.
\end{eqnarray*}

Thus, we have proved that, if
$\alpha= \epsilon_{i}$ and $m \geq 0$,
then  $ \langle \lambda_{S} + \rho,\tilde{\alpha}^{\vee}\rangle \notin
-{\Z}_{+}$.

Now, let's consider the case $\alpha= -\epsilon_{i}$ and $m \geq 1$.
Since $h_{i_{k}}<0$, using relation (\ref{2.6.2})
we get $( \mu _{S}, \epsilon_{i})=h_{i_{s}}+ \ldots \
+h_{i_{k}} \leq 0$.
It follows
\begin{eqnarray*}
&& \langle \lambda_{S} + \rho,\tilde{\alpha}^{\vee}\rangle =
2\left(m(l+\frac{1}{2})-( \bar{\rho}
, \epsilon_{i} )- ( \mu _{S} , \epsilon_{i} ) \right) \geq \\
&& \qquad \qquad \quad \ \ \geq 2l+1-(2l-2i+1)=2i>0,
\end{eqnarray*}
which implies that
$ \langle \lambda_{S} + \rho,\tilde{\alpha}^{\vee}\rangle \notin
-{\Z}_{+}$ holds in this case. \\
{\bf Case 3.}: The root $\alpha$ is of the form
$\alpha= \pm  (\epsilon_{i} + \epsilon_{j})$, for $i,j=1, \ldots ,l$,
$i<j$.

Then $( \bar{\rho}, \epsilon_{i} +
\epsilon_{j})=2l-j-i+1$. Let $s,t \in \{ 1, \ldots ,k \}$ be indices
such that $S \cap \{ i,i+1, \ldots ,j-1 \}= \{ i_{s}, \ldots ,i_{t}
\}$. Clearly, $i_{s} \geq i$ and $i_{t} \leq j-1$. Furthermore,
$( \mu _{S}, \epsilon_{i} + \epsilon_{j})=h_{i_{s}}+ \ldots \ +h_{i_{t}}+
2(h_{i_{t+1}}+ \ldots \ +h_{i_{k}})$.

First consider the case
$\alpha= \epsilon_{i}+ \epsilon_{j}$, $i<j$ and $m \geq 0$.

Then
\begin{eqnarray*}
&&\langle \lambda_{S} + \rho,\tilde{\alpha}^{\vee}\rangle =
m(l+\frac{1}{2})+( \bar{\rho}
, \epsilon_{i}+ \epsilon_{j} )+ ( \mu _{S} , \epsilon_{i}+
\epsilon_{j}) \\
&& \qquad \qquad \quad \ \ \geq 2l-j-i+1 + ( \mu _{S} , \epsilon_{i}+
\epsilon_{j}).
\end{eqnarray*}

Suppose that $t-s+1$ is even. If $k-t$ is also even,
then using relation (\ref{2.6.2}) we obtain
\begin{eqnarray*}
&& ( \mu _{S}, \epsilon_{i}+\epsilon_{j}) \geq -(i_{t}-i_{s})-
2(i_{k}-i_{t+1}),
\end{eqnarray*}
which implies
\begin{eqnarray*}
&&\langle \lambda_{S} + \rho,\tilde{\alpha}^{\vee}\rangle
\geq 2l-j-i+1 -(i_{t}-i_{s})- 2(i_{k}-i_{t+1})\\
&& \qquad \qquad \quad \ \ = 2(l-i_{k})+(i_{s}-i)+(i_{t+1}-i_{t})
+(i_{t+1}-j)+1>0.
\end{eqnarray*}
If $k-t$ is odd,
then using relation (\ref{2.6.2}) we obtain
\begin{eqnarray*}
&& ( \mu _{S}, \epsilon_{i}+\epsilon_{j}) \geq -(i_{t}-i_{s})-
2(i_{k-1}-i_{t+1})+2h_{i_{k}}, \\
&& \qquad \qquad \quad \ \ = -(i_{t}-i_{s})-
2(i_{k-1}-i_{t+1})+2(-l+i_{k}+\frac{1}{2}),
\end{eqnarray*}
which implies
\begin{eqnarray*}
&&\langle \lambda_{S} + \rho,\tilde{\alpha}^{\vee}\rangle
\geq 2l-j-i+1-(i_{t}-i_{s})-
2(i_{k-1}-i_{t+1})-2l+2i_{k}+1\\
&& \qquad \qquad \quad \ \ = 2(i_{k}-i_{k-1})+ (i_{t+1}-i_{t})
+(i_{t+1}-j)+(i_{s}-i)+2 >0.
\end{eqnarray*}

Suppose now that $t-s+1$ is odd. Then
\begin{eqnarray*}
( \mu _{S}, \epsilon_{i} + \epsilon_{j})=h_{i_{s}}+ \ldots \ +h_{i_{t-1}}+
(h_{i_{t}}+ 2h_{i_{t+1}}+ \ldots \ +2h_{i_{k}}).
\end{eqnarray*}
By using relation (\ref{2.6.2}) and system (\ref{2.6.1.1})
we get
\begin{eqnarray*}
&&( \mu _{S}, \epsilon_{i} + \epsilon_{j}) \geq -(i_{t-1}-i_{s})+
(h_{i_{t}}+ 2h_{i_{t+1}}+ \ldots \ +2h_{i_{k}}) \\
&& \qquad \qquad \quad \ \ =-(i_{t-1}-i_{s})+
(-l+i_{t}+\frac{1}{2}).
\end{eqnarray*}
It follows
\begin{eqnarray*}
&&\langle \lambda_{S} + \rho,\tilde{\alpha}^{\vee}\rangle
\geq 2l-j-i+1 -(i_{t-1}-i_{s})
-l+i_{t}+\frac{1}{2}\\
&& \qquad \qquad \quad \ \ = (l-j)+(i_{t}-i_{t-1})+
(i_{s}-i)+ \frac{3}{2}>0.
\end{eqnarray*}

Thus, we have proved that, if
$\alpha= \epsilon_{i} + \epsilon_{j}$, $i<j$ and $m \geq 0$,
then \linebreak $ \langle \lambda_{S} + \rho,\tilde{\alpha}^{\vee}\rangle \notin
-{\Z}_{+}$.

The only case left is
$\alpha= -(\epsilon_{i} + \epsilon_{j})$, $i<j$ and $m \geq 1$.
Then
\begin{eqnarray*}
&&\langle \lambda_{S} + \rho,\tilde{\alpha}^{\vee}\rangle =
m(l+\frac{1}{2})-( \bar{\rho}
, \epsilon_{i}+ \epsilon_{j} )- ( \mu _{S} , \epsilon_{i}+
\epsilon_{j}).
\end{eqnarray*}

If $t-s+1$ is even, then $( \mu _{S}, \epsilon_{i} +
\epsilon_{j}) \in \Z$, so if $m$ is odd,
then \linebreak $\langle \lambda_{S} + \rho,\tilde{\alpha}^{\vee}\rangle
\notin \Z$. Let $m$ be even. Then $m \geq 2$. Clearly
\begin{eqnarray*}
( \mu _{S}, \epsilon_{i} + \epsilon_{j})=(h_{i_{s}}+ \ldots \ +h_{i_{k}})+
(h_{i_{t+1}}+ \ldots \ +h_{i_{k}})\leq 0,
\end{eqnarray*}
which implies
\begin{eqnarray*}
&&\langle \lambda_{S} + \rho,\tilde{\alpha}^{\vee}\rangle \geq
2(l+\frac{1}{2})-( 2l-j-i+1 )=i+j>0.
\end{eqnarray*}

If $t-s+1$ is odd, then using
relation (\ref{2.6.2}) and system (\ref{2.6.1.1})
we obtain
\begin{eqnarray*}
( \mu _{S}, \epsilon_{i} + \epsilon_{j})=(h_{i_{s}}+ \ldots \ +h_{i_{t-1}})+
(h_{i_{t}}+2h_{i_{t+1}} \ldots \
+2h_{i_{k}}) \leq -l+i_{t}+\frac{1}{2}.
\end{eqnarray*}
It follows
\begin{eqnarray*}
&&\langle \lambda_{S} + \rho,\tilde{\alpha}^{\vee}\rangle \geq
(l+\frac{1}{2})-( 2l-j-i+1
)-(-l+i_{t}+\frac{1}{2})= \\
&& \qquad \qquad \quad \ \ =(j-1-i_{t})+i>0.
\end{eqnarray*}
We have proved that, if
$\alpha=- (\epsilon_{i} + \epsilon_{j})$, $i<j$ i $m \geq 1$,
then \linebreak $ \langle \lambda_{S} + \rho,\tilde{\alpha}^{\vee}\rangle \notin
-{\Z}_{+}$.

Thus, we have proved the relation (\ref{2.6.3.1}). Moreover, it can be easily checked
that coroots
\begin{eqnarray*}
&& (\delta - \alpha _{i_{j}})^{\vee}, \ j=1, \ldots , k \\
&& \alpha _{i_{j}}^{\vee}+\alpha _{i_{j}+1}^{\vee}+ \ldots + \alpha _{i_{j+1}}^{\vee},
\ j=1, \ldots , k-1 \\
&& \alpha _{i}^{\vee}, \ i \notin S, \ i \in \{1,2, \ldots ,l \}
\end{eqnarray*}
are elements of the set $\hat{\Delta}^{\vee \mbox{\scriptsize{re}}}_{\lambda _{S}}$
which implies
${\Q} \hat{\Delta}^{\vee \mbox{\scriptsize{re}}}_{\lambda _{S}}={\Q}
\hat{\Pi}^{\vee}$, and the relation (\ref{2.6.3.2}) is proved. $\;\;\;\;\Box$

\begin{thm} \label{t.2.6.5}
Let $M$ be a weak $L(-l+\frac{3}{2},0)$-module from the category
$\mathcal{O}$. Then $M$ is completely reducible.
\end{thm}
{\bf Proof:} Let $L(\lambda)$ be some irreducible subquotient
of $M$. Then $L(\lambda)$ is a \linebreak $L(-l+\frac{3}{2},0)$-module,
and Theorem \ref{c.2.6.3} implies that there exists \linebreak
$S \subseteq \{1,2, \ldots ,$ $ l-1 \}$ such that
$\lambda=(-l+\frac{3}{2}) \Lambda_{0}+ \mu _{S}$ or
$\lambda=(-l+\frac{3}{2}) \Lambda_{0}+ \mu _{S}'$.
It follows from Lemma \ref{l.2.6.4} that such $\lambda$ is admissible.
Proposition \ref{t.KW2} now implies that $M$ is completely reducible.
$\;\;\;\;\Box$

\bibliography{thesis}
\bibliographystyle{plain}

\vskip 1cm

Department of Mathematics, University of Zagreb, Bijeni\v{c}ka 30,
\linebreak 10000 Zagreb, Croatia

E-mail address: perse@math.hr

\end{document}